\documentclass[a4paper,12pt]{article}

\usepackage{amsmath,amstext,amscd,amsthm,amsfonts}
\usepackage{amssymb}
\newtheorem{theorem}{Theorem}[section]
\newtheorem{lemma}[theorem]{Lemma}

\newtheorem{definition}[theorem]{Definition}
\newtheorem{corollary}[theorem]{Corollary}
\theoremstyle{remark}
\newtheorem{rk}[theorem]{Remark}

\def\ad{\mathop{\rm ad}\nolimits}

\def\Aut{\mathop{\rm Aut}\nolimits}

\def\SL{\mathop{\rm SL}\nolimits}

\def\Id{\mathop{\rm Id}\nolimits}

\def\min{\mathop{\rm min}\nolimits}

\def\gr{\mathop{\rm gr}\nolimits}
\def\modulo{\mathop{\rm mod}\nolimits}
\def\Cas{\mathop{\rm Cas}\nolimits}

\title{Isomorphism Problems of Noncommutative Deformations of Type $D$ Kleinian singularities}

\author{Paul Levy //
        paul.levy@epfl.ch}

\begin{document}

\maketitle{}{}

\begin{abstract}
We construct all possible noncommutative deformations of a Kleinian singularity ${\mathbb C}^2/\Gamma$ of type $D_n$ in terms of generators and relations, and solve the problem of when two deformations are isomorphic.
We prove that all isomorphisms arise naturally from the action of the normalizer $N_{\SL(2)}(\Gamma)$ on ${\mathbb C}/\Gamma$.
We deduce that the moduli space of isomorphism classes of noncommutative deformations in type $D_n$ is isomorphic to a vector space of dimension $n$.
\end{abstract}

\section{Introduction}

Let $V$ be a complex vector space of dimension 2 and let $\Gamma$ be a finite subgroup of $\SL(V)$.
Such subgroups are classified: up to conjugacy, they are in one-to-one correspondence with the simply-laced Dynkin diagrams $A_n (n\geq 1),D_n (n\geq 4),E_6,E_7,E_8$.
Let $\Delta$ be the Dynkin diagram of $\Gamma$.
The quotient $V/ \Gamma$, which has coordinate ring ${\mathbb C}[V]^\Gamma$ and embeds as a hypersurface in ${\mathbb A}^3$ is a {\it Kleinian singularity} or {\it rational double point} of type $\Delta$.
The Dynkin diagram $\Delta$ arises as the exceptional configuration of the minimal resolution of the singularity $V/\Gamma$ (see \cite[\S 6]{slodowy}), or as the type of the McKay graph of $\Gamma$, which is isomorphic to the extended Dynkin diagram $\hat\Delta$ (\cite{mckay}).

It follows from the identification of $V/\Gamma$ with a hypersurface in ${\mathbb C}^3$ that there is a Poisson bracket on ${\mathbb C}[V]^\Gamma$, and an associated Poisson structure on the polynomial ring ${\mathbb C}[X,Y,Z]$.
In \cite{cbh}, Crawley-Boevey and Holland constructed a family of (in general non-commutative) deformations ${\cal O}_\lambda$ of (the Poisson bracket on) ${\mathbb C}[V]^\Gamma$, parametrised by $\lambda\in Z({\mathbb C}\Gamma)$.
This generalized work of Hodges \cite{hodges} and Smith \cite{smith} who constructed deformations of, respectively, a Kleinian singularity of type $A$ and the corresponding Poisson structure on ${\mathbb C}[X,Y,Z]$.
It is perhaps a little surprising that noone has attempted to describe the possible deformations of the non-type $A$ singularities in terms of generators and relations.
In Sect. 1 we carry this out for type $D$.
We show that construct all noncommutative deformations of a Kleinian singularity of type $D_n$, parametrised by a pair $(Q,\gamma)$ where $Q(t)$ is a monic polynomial of degree $(n-1)$ and $\gamma\in{\mathbb C}$.
We denote the corresponding associative algebras by $D(Q,\gamma)$ (Def. \ref{DQgamma}).
We also classify the noncommutative deformations of the corresponding Poisson structure on ${\mathbb C}[X,Y,Z]$, which are parametrised by $\gamma\in{\mathbb C}$ and a polynomial $P$ with leading term $(n-1)t^{n-2}$.
We denote the associative algebras thus produced by $H(P,\gamma)$.
(One could perform this process for the exceptional types.
The calculations are rather detailed, but not impossible.)

Let ${\mathfrak g}$ be a complex simple Lie algebra with Dynkin diagram $\Delta$ and let $E$ be a subregular nilpotent element of ${\mathfrak g}$.
Choose an $\mathfrak{sl}(2)$-triple $\{ H,E,F\}\subset{\mathfrak g}$ containing $E$.
It was proved by Brieskorn \cite{brieskorn} that the intersection of the {\it Slodowy slice} $E+{\mathfrak z}_{\mathfrak g}(F)$ with the nilpotent cone of ${\mathfrak g}$ is isomorphic to $V/ \Gamma$.
(This was generalized to the case of non-simply-laced $\Delta$ by Slodowy \cite{slodowy}.)
In fact, one can carry out the same process for an arbitrary nilpotent orbit.
In \cite{premslice}, Premet proved that all singularities constructed in this way have natural non-commutative deformations (see also \cite{gg}).
It was recently proved (by Arakawa \cite{arakawa} for regular $E$, by De Sole-Kac et al for general $E$ \cite[Appendix]{desole-kac}) that Premet's deformations are isomorphic to the {\it finite (quantum) $W$-algebras} of mathematical physics, constructed via quantum Hamiltonian reduction and the BRST cohomology (see de Boer and Tjin \cite{deboertjin}).

The associative algebras constructed by Hodges and Smith have straightforward presentations in terms of generators and relations: if $a$ and $P$ are polynomials then let $A(a)$ (resp. $R(P)$) be the algebra with generators $e,f,h$ (resp. $H,A,B$) and relations $he-eh=e,hf-fh=-f,ef=a(h-1),fe=a(h)$ (resp. $HA-AH=A,HB-BH=-B,AB-BA=P(H)$).
If $P(t)=a(t-1)-a(t)$, then clearly $A(a)$ is a quotient of $R(P)$.
The algebras of Smith are sometimes called {\it generalized Weyl algebras} because of their similarity to the first Weyl algebra.
The problem of when two algebras $A(a_1),A(a_2)$ are isomorphic was solved by Bavula and Jordan in \cite{bav-jord}: the isomorphisms are precisely the `obvious' ones.
Namely, $A(a_1)\cong A(a_2)$ if and only if $a_1(t)=\eta a_2(\tau\pm t)$ for some $\eta\in{\mathbb C}^\times,\tau\in{\mathbb C}$.
The analagous isomorphism theorem for the generalized Weyl algebras then follows.
In Sect. 2 we tackle the isomorphism problem for the algebras $D(Q,\gamma)$.
If $n>4$ then (assuming $Q$ monic) only the isomorphisms $D(Q,\gamma)\cong D(Q,-\gamma)$ occur (Thm. \ref{geq4}).
It follows that the moduli space of isomorphism classes of deformations of (the Poisson bracket on) a Kleinian singularity of type $D_n$ is isomorphic to affine $n$-space, generalizing the result for type $A$ (Cor. \ref{moduli}).
In Sect. 3 we carry out specific calculations for the case $n=4$, where the situation is rather more interesting.
Here there are six sets of isomorphisms, corresponding to the six elements of $N_{\SL(V)}(\Gamma)/\Gamma\cong S_3$ (Thm. \ref{main}).
The moduli space of isomorphism classes is isomorphic to a vector space of dimension 4 (Cor. \ref{moduli4}).
We also apply our results on $D(Q,\gamma)$ to solve the problem of when two algebras $H(P,\gamma),H(\tilde{P},\tilde\gamma)$ are isomorphic (Thm. \ref{Hn4} and Thm. \ref{H4}).

Our methods share a certain similarity with those of Bavula and Jordan \cite{bav-jord}, who adapted Dixmier's approach to solving the isomorphism problem for the case $\deg a=1,2$ (\cite{dixmier,dixmier2}).
In particular, we construct a filtration of $D(Q,\gamma)$ by the additive monoid ${\mathbb Z}_{\geq 0}\times{\mathbb Z}_{\geq 0}$ (which compares to Bavula and Jordan's family of filtrations of $A(a)$ \cite[Thm. 3.14]{bav-jord}) and exploit a peculiar property of one of the generators for $D(Q,\gamma)$ to analyse its possible images in the corresponding graded algebra (Lemma \ref{reducetov}).
On the other hand, at this point we diverge sharply from the path of \cite{bav-jord}, since elimination of the remaining cases requires an in-depth study of certain expressions involving commutators.
However, one advantage of this analysis is the explicit construction of the `non-trivial' isomorphisms in type $D_4$ (Definition \ref{isos}).

Unfortunately, we have as yet been unable to determine whether every algebra $D(Q,\gamma)$ is isomorphic to some ${\cal O}_\lambda$.
Given our description of the isomorphisms as essentially arising from the normalizer of $\Gamma$ in $\SL(V)$ (which can be identified with the group of graph automorphisms of a root system of type $D_n$) we expect this to be the case.

{\it Notation.}
We denote by $[x,y]$ for the commutator product $xy-yx$.
If $m$ and $j$ are positive integers, then $[m/j]$ will denote the integer part of $m/j$.

{\it Acknowledgement.}
Much of the research for this paper was carried out while at Aarhus University under a fellowship from the ``Lie GRITS" European Union research network.
I am grateful both to Lie GRITS and to my Danish hosts for this opportunity.
The article was completed at the \'Ecole Polytechnique F\'ed\'erale de Lausanne.
I would like to thank Alexander Premet for alerting me to the problem and for several helpful conversations.
I would also like to express my appreciation of the advice of Jens Carsten Jantzen and Iain Gordon.

\section{Generators and Relations}

Let $V$ be a complex vector space of dimension 2.
Identify $\SL(V)$ with $\SL(2)$ by choice of a basis for $V$, and let $x,y$ be the corresponding coordinate functions on $V$.
Up to conjugacy, there is a unique {\it binary dihedral group} $\Gamma\subset\SL(V)$ of order $4(n-1)$ for each $n\geq 3$.
Following \cite{slodowy}, we choose the following generators for $\Gamma$:
$$\sigma=\left(\begin{array}{cc}
\zeta & 0 \\
0 & \zeta^{-1}
\end{array}\right)\;\mbox{and}\;\tau=\left(\begin{array}{cc}
0 & 1 \\
-1 & 0
\end{array}\right),\;\;\mbox{where }\zeta=e^{\pi i/(n-1)}.$$

The quotient $V/ \Gamma$ is a Kleinian singularity or rational double point of type $D_{n+1}$.
It is easy to see that the coordinate ring ${\mathbb  C}[V]^\Gamma$ is generated by $x^2y^2,(x^{2(n-1)}+y^{2(n-1)})$ and $xy(x^{2(n-1)}-y^{2(n-1)})$, hence is isomorphic to ${\mathbb C}[X,Y,Z]/(X^n+XY^2+Z^2)$.

Recall that a Poisson algebra is a commutative algebra $B$ endowed with a {\it Poisson bracket} $\{ .\, ,.\}$ satisfying:

(i) $(B,\{ .\, ,.\})$ is a Lie algebra,
 
(ii) $\{b,.\}$ and $\{.\, ,b\}$ are derivations of $B$ for each $b\in B$.

Any polynomial $\phi\in{\mathbb C}[X,Y,Z]$ induces a Poisson algebra structure on ${\mathbb C}[X,Y,Z]$, which we denote $\{ .\, ,.\}_\phi$, such that:
$$\{X,Y\}_\phi=\partial\phi/ \partial Z, \{X,Z\}_\phi=-\partial\phi/ \partial Y, \{Y,Z\}_\phi=\partial\phi/ \partial X.$$

Moreover, since $(\phi)\subset{\mathbb C}[X,Y,Z]$ is a Poisson ideal, there is an induced Poisson bracket on the quotient ${\mathbb C}[X,Y,Z]/(\phi)$.
In the case $\phi= X^n+XY^2+Z^2$, we will first construct (all possible) non-commutative deformations of the Poisson bracket on ${\mathbb C}[X,Y,Z]$.
We denote the algebras thus produced by $H(P,\gamma)$, parametrised by a polynomial $P$ of the form $nt^{n-1}+\ldots$ and a scalar $\gamma$.
The $H(P,\gamma)$ are the type $D$ analogues of the generalized Weyl algebras constructed by Smith \cite{smith}.
In Lemma \ref{centreisp} we show that the centre of $H(P,\gamma)$ is a polynomial ring ${\mathbb C}[\Omega]$ on one generator, and provide a precise description of $\Omega$.
The various factor algebras $H(P,\gamma)/(\Omega-c)$ thus determine all possible non-commutative deformations of the Kleinian singularity ${\mathbb C}[X,Y,Z]/(\phi)$ of type $D_{n+1}$.

For our purposes, a non-commutative deformation of a Poisson algebra $(A_0,\{ .\, ,.\})$ is an associative ${\mathbb C}[[t]]$-algebra ${\cal A}$, free as a ${\mathbb C}[[t]]$-module, such that:

(a) There is an isomorphism $\pi:{\cal A}/t{\cal A}\longrightarrow A_0$ of associative (commutative) algebras,

(b) For any $x,y\in {\cal A}$, $\pi(xy-yx+t{\cal A})=\{\pi(x+t{\cal A}),\pi(y+t{\cal A})\}$.

Note that freeness implies that any lift of a basis of $A_0$ is a (${\mathbb C}[[t]]$-)basis of ${\cal A}$.
But it follows that if $A_0$ is a Kleinian singularity of type $D_{n+1}$ (resp. the corresponding Poisson algebra on ${\mathbb C}[X,Y,Z]$) then any deformation ${\cal A}$ of $A_0$ possesses a set $U,V,W$ of generators such that $\{ U^iV^jW^\epsilon \,:\, i,j\geq 0,\epsilon\in\{ 0,1\}\,\}$ (resp. $\{ U^iV^jW^k\,:\, i,j,k\geq 0\}$) is a ${\mathbb C}[[t]]$-basis for ${\cal A}$.
Moreover, for any $\alpha,\beta\in{\mathbb C}^\times$ the quotients ${\cal A}/ (t-\alpha){\cal A}$ and ${\cal A}/ (t-\beta){\cal A}$ are naturally isomorphic, by appropriate scaling of the images of $U,V,W$.
Hence we can (and will) abuse terminology and refer to the quotient $A={\cal A}/(t-1){\cal A}$ as the  deformation of $A_0$.
In less formal language, a noncommutative deformation of $A_0$ is a filtered associative algebra $A$ satisfying the appropriate condition as above on a generating set, such that $\gr A=A_0$ and $\gr [x,y]=\{ \gr x,\gr y\}$.

For the moment we wish to determine all noncommutative deformations of the Poisson algebra $({\mathbb C}[X,Y,Z],\{ .\, ,.\}_\phi)$.
Hence suppose $A$ has generators $U,V,W$ such that $U^iV^jW^k$ is a basis.
We require that $\gr U=X$, $\gr V=Y$, and $\gr W=Z$: hence the filtration on $A$ satisfies $U\in A_4\setminus A_2,V\in A_{2n-2}\setminus A_{2n-4}$, and $W\in A_{2n}\setminus A_{2n-2}$.
Moreover, with respect to this filtration $[U,V]=2W+\mbox{lower terms}$, $[U,W]=-2UV+\mbox{lower terms}$, and $[V,W]=V^2+nU^{n-1}+\mbox{lower terms}$.
We wish to find the possible expressions for these commutators satisfying the Jacobi identity.
But we may clearly replace $U$ (resp. $V,W$) by an equivalent element modulo the scalars (resp. $A_{2n-4},A_{2n-2}$).

Hence, after substituting for $W$, we assume that $[U,V]=2W$.
Now $[U,W]=-2UV+\alpha W+\beta V+p(U)$ for some polynomial $p\in{\mathbb C}[t]$ of degree $\leq (n+1)/2$.
Substituting $(U-\beta/2)$ for $U$, we may assume that $\beta=0$.
Let $p=tq+\gamma$ for $q\in{\mathbb C}[t],\gamma\in{\mathbb C}$.
Replacing $V$ by $(V-q(U)/2)$, we may assume that $[U,W]=-2UV+\alpha W+\gamma$.
Finally, there exist polynomials $P,m_1,m_2\in{\mathbb C}[t]$ with $m_1$ of degree $\leq (n-3)/2$, $m_2$ of degree $\leq (n-2)/2$ and $P$ with leading term $nt^{n-1}$ such that $[V,W]=V^2+P(U)+m_1(U)W+m_2(U)V$.
But the Jacobi identity now requires that $[U,[V,W]]=[V,[U,W]]$, hence that $m_1=m_2=0$ and $\alpha=2$.

\begin{definition}
Let $P(t)$ be a polynomial of degree $(n-1)$ and let $\gamma\in{\mathbb C}$.
The algebra $H(P,\gamma)$ has generators $U,V,W$ and relations $[U,V]=2W$, $[U,W]=-2UV+2W+\gamma$ and $[V,W]=V^2+P(U)$.
\end{definition}

This definition does not require $P$ to have leading term $nt^{n-1}$; but by scaling generators $(U,V,W)\mapsto (U,\alpha V,\alpha W)$ we can easily see that $H(P,\gamma)$ is isomorphic to $H(\alpha^2 P,\alpha\gamma)$.
It follows immediately from the above discussion:

\begin{lemma}
Let $A$ be a noncommutative deformation of $({\mathbb C}[X,Y,Z],\{.\, ,.\}_\phi)$, where $\phi=X^n+XY^2+Z^2$.
Then $A$ is isomorphic (as a filtered algebra) to $H(P,\gamma)$ for some polynomial $P(t)=nt^{n-1}+\ldots$ and $\gamma\in{\mathbb C}$.
\end{lemma}

We now describe the centre of $H(P,\gamma)$.
To do this we need a little preparation.
By definition $[U,V]=2W$ and $[U,W]=-2UV+2W+\gamma$.
It follows that there exist polynomials $\alpha_n,\beta_n\in{\mathbb C}[t]$ such that $[U^n,V]=\alpha_n(U)[U,V]+\beta_n(U)[U,W]$.
Indeed, then $$\begin{array}{rl} [U^{n+1},V]= & U^n[U,V]+(\alpha_n(U)[U,V]+\beta_n(U)[U,W])U \\
= & (U^n+\alpha_n(U)+2U\beta_n(U))[U,V]+((U-2)\beta_n(U)-2\alpha_n(U))[U,W]\end{array}$$ 
Hence $\alpha_{n+1}=t^n+t\alpha_n+2t\beta_n$ and $\beta_{n+1}=(t-2)\beta_n-2\alpha_n$.
To solve these difference equations, let $\iota:{\mathbb C}[t]\hookrightarrow{\mathbb C}[s]$ be the algebra embedding which sends $t$ to $-s(s+1)$.
For $f\in{\mathbb C}[t]$ (temporarily) denote by $\overline{f}$ the image $\iota(f)$.
Let $\rho_n=\overline{\alpha_n}-s\overline{\beta_n}$ and let $\mu_n=\overline{\alpha_n}+(s+1)\overline{\beta_n}$.
A straightforward calculation shows that $\rho_{n+1}=(-s(s+1))^n-s(s-1)\rho_n$ and $\mu_{n+1}=(-s(s+1))^n-(s+1)(s+2)\mu_n$.
But $\rho_1=\mu_1=1$, hence $\rho_n=((-s(s-1))^n-(-s(s+1))^n)/2s$ and $\mu_n=((-s(s+1))^n-(-(s+1)(s+2))^n)/2(s+1)$.
It would be straightforward to write down explicit expressions for $\alpha_n$ and $\beta_n$, but this will suffice for our purposes.
Let $\rho,\mu:{\mathbb C}[s]\rightarrow{\mathbb C}[s]$ be the linear maps given by $\rho(p)=(p(-s)-p(s))/2s$ and $\mu(p)=(p(-(s+1))-p(s+1))/2(s+1)$.
We note that for any $f\in{\mathbb C}[t]$ there exist unique polynomials $\alpha(f),\beta(f)\in{\mathbb C}[t]$ such that $\overline{\alpha(f)}-s\overline{\beta(f)}=\rho(\overline{f})$.
Indeed, by the above discussion there exist unique $\alpha(f),\beta(f)$ such that $\overline{\alpha(f)}-s\overline{\beta(f)}=\rho(\overline{f})$ and $\overline{\alpha(f)}+(s+1)\overline{\beta(f)}=\mu(\overline{f})$.
But $-s(s+1)$ is stable under the algebra endomorphism of ${\mathbb C}[s]$ which sends $s$ to $-(s+1)$, hence the first condition implies the second.
Hence we introduce the linear maps $\alpha,\beta:{\mathbb C}[t]\rightarrow{\mathbb C}[t]$ such that $\overline{\alpha(f)}-s\overline{\beta(f)}=\rho(\overline{f})$ for all $f\in{\mathbb C}[t]$.

\begin{lemma}\label{ucomms}
(a) $[f(U),V]=\alpha(f)(U)[U,V]+\beta(f)(U)[U,W]$.

(b) $[f(U),W]=-U\beta(f)(U)[U,V]+(\alpha(f)+\beta(f))(U)[U,W]$.
\end{lemma}

\begin{proof}
The first assertion is an immediate consequence of the discussion in the paragraph above.
For (b), we note that $2[U,W]=[U,[U,V]]=2\alpha(f)(U)[U,W]-2U\beta(f)(U)[U,V]+2\beta(f)(U)[U,W]$.
\end{proof}

\begin{lemma}\label{centreisp}
Let $Q$ be a monic polynomial, unique up to addition of scalars, such that $Q(-s(s-1))-Q(-s(s+1))=(s-1)P(-s(s-1))+(s+1)P(-s(s+1))$ and let $\Omega=Q(U)+UV^2+W^2-2WV-\gamma V\in H=H(P,\gamma)$.
Then $Z(H)={\mathbb C}[\Omega]$.
\end{lemma}

\begin{proof}
Let $h$ be an element of $H$ of the form $Q(U)+UV^2+W^2+\alpha WV+\beta V^2+p_1(U)W+p_2(U)V$, where $Q\in{\mathbb C}[t]$ is monic of degree $n$ and $p_1,p_2\in{\mathbb C}[t]$ are polynomials of degrees $\leq (n-1)/2$ and $\leq n/2$ respectively.
To find the possible $Q,p_1,p_2$ such that $h\in Z(H)$ we have only to find the conditions under which $[z,U]=[z,V]=0$ (since then $[z,[U,V]]=0$, hence $z\in Z(H)$).
It is easy to see that $$[U,UV^2+W^2]=2UWV-2WUVv+4W^2+2\gamma W=[U,2WV+\gamma V]$$
It follows that $[U,h]=0$ if (and only if, though this is unnecessary) $h=Q(U)+UV^2+W^2-2WV-\gamma V$, for some monic polynomial $Q$ of degree $n$.
Assume $h$ is of this form.
To determine when $[h,V]=0$ we apply Lemma \ref{ucomms}.
By a straightforward calculation, $$\begin{array}{rl}[UV^2+W^2-2WV,V]= & -[V,[V,W]]-2P(U)W+[P(U),W] \\ = & [P(U),V]+[P(U),W]-P(U)[U,V]\end{array}$$
But therefore $h\in Z(H)$ if and only if $[Q(U)+P(U),V]+[P(U),W]=P(U)[U,V]$.
By Lemma \ref{ucomms}, $\alpha(Q)=P+t\beta(P)-\alpha(P)$ and $\beta(Q)=-(\alpha +2\beta)(P)$.
It follows that $\overline{\alpha(Q)}-s\overline{\beta(Q)}=\overline{P}+(s-1)(\alpha-s\beta)(\overline{P})$.
Hence $Q$ is the unique polynomial modulo addition of scalars such that $\overline{Q}(-s)-\overline{Q}(s)=(s-1)\overline{P}(-s)+(s+1)\overline{P}(s)$.

This proves that $\Omega=Q(U)+UV^2+W^2-2WV-\gamma V\in Z(H)$.
Let $B$ be the Poisson algebra $({\mathbb C}[X,Y,Z],\{.\, ,.\}_\phi)$.
It is well-known (and easy to check) that the Casimir elements $\Cas B=\{ f\in B\,:\, \{f,g\}=0\,\forall g\in B\}={\mathbb C}[\phi]$.
It is easy to see that if $h\in Z(H)$ then $\gr h\in\Cas B$.
Suppose therefore that $h\in Z(H)$, but $h\not\in {\mathbb C}[\Omega]$.
We may assume that the degree of $h$ is minimal subject to this condition.
Then $\gr h\in Cas B$, hence $\gr h=\xi\phi^i$ for some $i$ and some $\xi\in{\mathbb C}^\times$.
But now $h-\xi\Omega^i\in Z(H)$ has degree strictly less than $h$, which contradicts our original assumption.
\end{proof}

We note that the condition on $Q,P$ is equivalent to the condition: \begin{eqnarray}\label{QP} Q(-s(s+1))+(s+1)P(-s(s+1)) \mbox{ is an even polynomial in $s$}\end{eqnarray}
Moreover, for each monic polynomial $Q(t)$ there is a unique $P(t)$ satisfying (\ref{QP}), necessarily with leading term $nt^{n-1}$ (where $n$ is the degree of $Q$).

\begin{definition}\label{DQgamma}
Let $Q(t)$ be a polynomial of degree $n$ and let $\gamma\in{\mathbb C}$.
We define $D(Q,\gamma)$ to be the associative algebra with generators $u,v,w$ and relations:
$$[u,v]=2w,\;\;[u,w]=-2uv+2w+\gamma,\;\;[v,w]=v^2+P(u)\;\;\mbox{and}\;\;Q(u)+uv^2+w^2-2wv-\gamma v=0$$

where $P(t)$ is the unique polynomial of degree $(n-1)$ such that $$Q(-s(s-1))-Q(-s(s+1))=(s-1)P(-s(s-1))+(s+1)P(-s(s+1)).$$
\end{definition}

In common with the convention for type $A$, we have not assumed that $Q$ is monic in the above definition.
But the change of generators $(u,v,w)\mapsto (u,\xi v,\xi w)$ gives a natural isomorphism $D(Q,\gamma)\cong D(\xi^2 Q,\xi\gamma)$.
Hence any such algebra $D(Q,\gamma)$ is isomorphic to some $D(Q_0,\gamma_0)$ with $Q_0$ monic.

\section{The Isomorphism Problem}

Recall that if $A$ is any ${\mathbb Z}$-filtered algebra, then there is a uniquely defined degree function on non-zero elements of $A$: $\deg x=\min_{x\in A_i}i$.
Fix a monic polynomial $Q(t)$ of degree $n\geq 3$ and $\gamma\in{\mathbb C}$, and let $A=D(Q,\gamma)$.
Let $P(t)$ be the unique polynomial such that $Q(-s(s+1))+(s+1)P(-s(s+1))$ is even in $s$.
By construction $A$ is a ${\mathbb Z}$-filtered algebra such that $u$ has degree 4, $v$ has degree $2 n-2$ and $w$ has degree $2n$.
Specifically, $\{ u^iv^jw^\epsilon\,:\, i,j\in{\mathbb Z}\geq 0,\epsilon\in\{ 0,1\}\}$ is a basis for $A$ and $\deg \sum_{i,j,\epsilon}a_{ij\epsilon}u^iv^jw^\epsilon=\max_{a_{ij\epsilon}\neq 0}(4i+(2n-2)j+2n\epsilon)$.

However, for any $N>n$ we can also define a filtration on $A$ with degree function $\deg\sum_{i,j,\epsilon}a_{ij\epsilon}u^iv^jw^\epsilon=\max_{a_{ij\epsilon}\neq 0}(4i+(2N-2)j+2N\epsilon)$.
To see this we have only to check that if $x,y\in A$ then $\deg xy\leq \deg x+\deg y$.
Hence suppose $x=\sum a_{ij\epsilon}u^iv^jw^\epsilon$ and $y=\sum b_{ij\epsilon}u^iv^jw^\epsilon$.
Then $4i+(2N-2)j+2N\epsilon\leq\deg x$ for all $a_{ij\epsilon}\neq 0$, and $4k+(2N-2)l+2N\eta\leq\deg y$ for all $b_{kl\eta}\neq 0$.
It follows that $4(i+k)+(2N-2)(j+l)+2N(\epsilon+\eta)\leq\deg x+\deg y$ for all $a_{ij\epsilon}b_{kl\eta}\neq 0$.
But $xy=\sum a_{ij\epsilon}b_{kl\eta}(u^{i+k}v^{j+l}w^{\epsilon+\eta}-u^{i+k}v^j[v^{l},w^\epsilon]w^\eta-u^i[u^k,v^jw^\epsilon]v^lw^\eta)$.
Hence by induction on $\deg x,\deg y$ we have only to show that the commutator relations $[u,v]=2w$, $[u,w]=-2uv+2w+\gamma$, $[v,w]=v^2+P(u)$ and the substitution $w^2=-Q(u)-uv^2-2vw+2v^2+2P(u)+\gamma v$ are of non-positive degree, that is, the terms on the right are of equal or lower degree than each of the terms on the left.
This is easily checked.
It will be extremely useful to us to consider the `limit as $N$ tends to infinity' of these filtrations.
Hence consider the additive monoid of pairs $(a,b)$ of non-negative integers, with the lex ordering $(a,b)>(a',b')$ if and only if $a>a'$ or $a=a'$ and $b>b'$.
Let $A_0^{(0)}={\mathbb C}\subset A$, let $A_a^{(b)}$ be the subspace of $A$ spanned by all monomials of the form $u^iv^jw^\epsilon$ with $(j+\epsilon,2i+\epsilon)\leq (a,b)$ and let $A_a^{(\infty)}=\cup_{b\geq 0} A_a^{(b)}$.
(We assume that $\epsilon\in\{ 0,1\}$, although this isn't strictly necessary.)
It is straightforward to check that the commutation relations $[u,v]=2w$, $[u,w]=-2uv+2w+\gamma$ and $[v,w]=v^2+P(u)$ satisfy $\deg [x,y]=\deg x+\deg y-(0,1)$ and that the equality $w^2=-uv^2-Q(u)+2wv+\gamma v$ replaces $w^2$ be a term of equal degree $(2,2)$ (congruent to $-uv^2$ modulo $A_2^{(1)}$).
It follows by the argument above that $A=\cup_{a,b\geq 0} A_a^{(b)}$ is a well-defined filtration of $A$, which we call the {\it limit filtration}.
It is easy to see that the corresponding graded algebra is isomorphic to ${\mathbb C}[X,Y,Z]/(XY^2+Z^2)$, where $X$ has degree $(0,2)$, $Y$ has degree $(1,0)$ and $Z$ has degree $(1,1)$.

Until further notice we fix the limit filtration on $A$.
It turns out to be significantly easier for us to calculate using the monomials $u^iwv^{j-1}$ rather than $u^iv^{j-1}w$.
(This does not effect our definition of the filtration since $wv=vw+$terms of lower degree.)
Hence we express elements of $A$ in terms of the basis $\{ u^iw^\epsilon v^j\,:\, i,j\geq 0,\epsilon=0,1\}$.
Since any subset of the ordered set ${\mathbb Z}_{\geq 0}\times {\mathbb Z}_{\geq 0}$ has a minimal element, there is a well-defined degree function on non-zero elements of $A$.
It is easy to see moreover that $u^iw^\epsilon v^j$ has degree $(a,2b)$ if and only if $i=b$, $j=a$ and $\epsilon=0$, and has degree $(a,2b+1)$ if and only if $a>0$, $\epsilon=1$, $i=b$ and $j=a-1$.
Hence each summand in the grading of $\gr A$ is of dimension 1.
We will write $x=\xi u^iw^\epsilon v^j+$lower terms to mean that $x$ is congruent to $\xi u^iw^\epsilon v^j$ modulo $\cup_{(a,b)<(j+\epsilon,2i+\epsilon)}A_a^{(b)}$ (implicitly assuming $\xi\neq 0$).
We refer to $\xi u^iw^\epsilon v^j$ as the `leading term' in $x$.

%
%
Note that $[u,v^m]=2mwv^{m-1}+$lower terms and $[u,wv^{m-1}]=-2muv^m+$lower terms, thus the cosets of $(\ad u)^j(v^m)$, $j\geq 0$ form a basis for $A_m^{(\infty)}/A_{m-1}^{(\infty)}$.
Define polynomials $F_m\in{\mathbb C}[S,T],\; m\in{\mathbb Z}_{\geq 0}$ by: $F_0=S$ and $F_m=(S^2-2m^2 S+m^2(m^2-1) + 4m^2T)$ for $m\geq 1$.
Clearly $\ad u$ and left multiplication by $u$, denoted $l_u$, commute.

\begin{lemma}\label{Fprod}
Let $x\in A$.
Then there exists $m$ such that $\prod_{i=0}^m F_i(\ad u,l_u)(x)=0$.
\end{lemma}

\begin{proof}
Since $\ad u$ and $l_u$ preserve each of the subspaces $A_m^{(\infty)}$, it is enough to show that $F_m(\ad u,l_u)(x)\in A_{m-1}^{(\infty)}$ for any $x\in A_m^{(\infty)}$.
But $A_m^{(\infty)}/A_{m-1}^{(\infty)}$ is spanned by (the cosets of) $(\ad u)^i(v^m)$, $i\in{\mathbb Z}_{\geq 0}$, hence it will suffice to show that $F_m(\ad u,l_u)(v^m)\in A_{m-1}^{(\infty)}$.

Clearly $$\begin{array}{rl}[u,v^m]= & 2wv^{m-1}+2vwv^{m-2}+\ldots+2v^{m-1}w \\ = & 2mwv^{m-1}+\sum_{j=1}^{m-1}[v^j,w]v^{m-1-j}\end{array}$$
But since $[v,w]\equiv v^2\;(\modulo A_0^{(\infty)})$, we deduce that $[v^j,w]v^{m-1-j}\equiv jv^m\;(\modulo A_{m-2}^{(\infty)})$.
It follows that $[u,v^m]\equiv 2mwv^{m-1}+m(m-1)v^m\;(\modulo A_{m-2}^{(\infty)})$.
Now $[u,wv^{m-1}]=-2uv^m+2wv^{m-1}+\gamma v^{m-1}+w[u,v^{m-1}]$.
Thus $[u,wv^{m-1}]\equiv -2uv^m+((m-1)(m-2)+2)wv^{m-1}+\gamma v^{m-1}+2(m-1)w^2v^{m-2}\;(\modulo A_{m-2}^{(\infty)})$.
But $w^2\equiv -uv^2+2wv+\gamma v\;(\modulo A_0^{(\infty)})$ by the defining relations for $A$.
Hence $[u,wv^{m-1}]\equiv -2muv^m+m(m+1)wv^{m-1}+(2m-1)\gamma v^{m-1}\;(\modulo A_{m-2}^{(\infty)})$.

This proves the statement about $v^m$: in fact we have shown that $F_m(\ad u,l_u)(v^m)\equiv 2m(2m-1)\gamma v^{m-1}\;(\modulo A_{m-2}^{(\infty)})$.
\end{proof}

\begin{lemma}\label{prodform}
Let $P(S,T)=\prod_{i=0}^m F_i(S,T)$.
If $P$ is written in the form $\sum_{i=0}^{2m+1} a_i(T)S^i$, then $a_{2m+1}=1$ and $\deg a_{2m+1-i}\leq i/2$.
\end{lemma}

\begin{proof}
Let ${\cal S}$ be the set of all polynomials in ${\mathbb C}[S,T]$ of the form $\sum_{i=0}^N a_i(T)S^i$, where $a_N\neq 0$ and $\deg a_{N-i}\leq i/2$.
The product of any two elements of ${\cal S}$ is also in ${\cal S}$, since the coefficient of $S^i$ in $(\sum a_j(T)S^j)(\sum b_l(T)S^l)$ is $\sum_{j+l=i}a_j(T)b_l(T)$.
But clearly $F_i\in{\cal S}$ for all $i$, hence $P\in{\cal S}$.
\end{proof}

Denote by $\gr_{lim} A$ the graded algebra of $A$ corresponding to the limit filtration, identified with ${\mathbb C}[X,Y,Z]/(XY^2+Z^2)$.
The Poisson bracket on $\gr_{lim}A$ satisfies $\{ X,Y\}=2Z$, $\{ X,Z\}=-2XY$, $\{ Y,Z\}=Y^2$.

\begin{lemma}\label{monomials}
Let $x$ be a monomial in $X,Y,Z$.
Then unless $x=Y^b$ or $x=X^b$, there exists some monomial $y$ in $X,Y,Z$ such that $\{ x,.\}^M(y)\neq 0$ for all $M\geq 0$.
\end{lemma}

\begin{proof}
Since $Z^2=-XY^2$, we have only to prove the lemma in the case where $x=X^aY^bZ^\epsilon$ with $\epsilon\in\{ 0,1\}$.
Suppose $x'=x^r$ for some $r\geq 2$.
Then $\{ x',.\}(y)=rx^{r-1}\{ x,y\}$.
Hence $\{ x',.\}^M(y)=r^Mx^{M(r-1)}\{ x,.\}^M(y)$.
It follows that we need only prove the lemma in the case where $x$ cannot be expressed as a power of any other monomial.
Note that if $(i,j)$ is the degree of $x$ in $\gr_{lim} A$ then this holds if and only if $i$ and $j$ are coprime.

Suppose first of all that $x=X^aY^b$ such that $b$ and $2a$ are coprime (and $ab\neq 0$).
By calculation $\{ X^aY^b,X^cY^d\}=2(ad-bc)X^{a+c-1}Y^{b+d-1}Z$.
Moreover, $\{ X^aY^b,X^cY^{d-1}Z\}=(b(2c+1)-2ad))X^{a+c}Y^{b+d}$.
Hence by our condition on $x$, $\{ x,y\}=0$ if and only if $y=x^r$ for some $r$.
We claim that $\{ x,.\}^{2i+1}(X^cY^d)=0$ if and only if $(c,d)=(ka+j,kb)$ for some $k\geq 0$, $0\leq j\leq i$ and $\{ x,.\}^{2i}(X^cY^{d-1}Z)=0$ if and only if $(c,d)=(ka+j,kb)$ for some $k\geq 0$, $0\leq j\leq (i-1)$.
This is true for $i=0$ by the above calculations.
Hence suppose we know our claim to be true for $(i-1)$.
Then $\{ x,.\}^{2i}(X^cY^{d-1}Z)=0$ if and only if $(b(2c+1)-2ad)\{ x,.\}^{2i-1}(X^{a+c}Y^{b+d})=0$.
By the induction hypothesis, this is true if and only if $(c,d)=(ka+j,kb)$ for some $k\geq 0$ and some $j$, $0\leq j\leq (i-1)$.
This proves the induction step for $X^cY^{d-1}Z$.
But now $\{ x,.\}^{2i+1}(X^cY^d)=0$ if and only if $2(ad-bc)\{ x,.\}^{2i}(X^{a+c-1}Y^{b+d-1}Z)=0$.
It follows from the above step that $\{ x,.\}^{2i+1}(X^cY^d)=0$ if and only if $(c,d)=k(a,b)$ or $(c-1,d)=(ka+j,kb)$ for some $j$, $0\leq j\leq i-1$.
This proves our claim.
Thus $\{ x,.\}^M(y)=0$ for some $M$ if and only if $y=X^{ka+c}Y^{kd}$ or $y=X^{kc+a}Y^{kd-1}Z$ for some $k\in{\mathbb N}$.
Therefore (for example) $\{ x,.\}^M(Y)\neq 0$ for all $M\geq 0$.

Suppose now that $x=X^aY^{b-1} Z$ where $b$ and $(2a+1)$ are coprime.
By the above $\{ x,X^cY^d\}=(2bc-d(2a+1))X^{a+c}Y^{b+d}$.
It follows that $\{ x,X^cY^{d-1}Z\}=(b(2c+1)-d(2a+1))X^{a+c}Y^{b+d-1}Z$.
Thus once more $\{ x,y\}=0$ if and only if $y=x^r$ for some $r$.
If $a$ and $b$ are not both zero, then it follows that $\{ x,.\}^M(Z)$ is a non-zero multiple of $X^{Ma}Y^{Mb+M}Z$ for each $M\geq 0$.
On the other hand, if $x=Z$ then $\{ x,.\}^M(Y)$ is a non-zero multiple of $Y^{M+1}$ for each $M\geq 1$.
This completes the proof.
\end{proof}

\begin{rk}
We note that it follows from the proof of Lemma \ref{monomials}, if $x$ and $y$ are monomials of coprime degrees $(m,i)$ and $(m',i')$, then $\{ x,y\}=\pm (mi'-im') z$, where $z$ is a monomial of degree $(m+m',i+i'-1)$.
\end{rk}

\begin{lemma}\label{reducetov}
Suppose $f$ is an element of $A$ satisfying the condition that for any $a\in A$ there exists $m$ such that $\prod_{i=0}^m F_i(\ad f,l_f)(a)=0$.
Then either $f\in{\mathbb C}[u]$ or there exist $r$ and $\xi\neq 0$ such that $f=\xi v^r+$lower terms.
\end{lemma}

\begin{proof}
Let $f$ be such an element, let $a\in A$ and suppose $\prod_0^m F_i(\ad f,l_f)(a)=0$.
Recall by Lemma \ref{prodform} that $P(S,T)=\prod_{i=0}^m F_i(S,T)$ is of the form $S^{2m+1}+a_{2m}S^{2m}+a_{2m-1}(T)S^{2m-1}+\ldots +a_0(T)$, where $\deg a_{2m+1-i}\leq i/2$.
Suppose $\gr f=x$ is not of the form $\xi X^i$ or $\xi Y^i$.
Then by Lemma \ref{monomials} there exists $y\in\gr_{lim} A$ such that $\{x,.\}^M(y)\neq 0$ for all $M\geq 0$.
Let $a\in A$ be such that $\gr a=y$.
Then it is easy to see that $\deg(\ad f)^m(y)=\deg \{x,.\}^m(y)$.
Let $\deg\gr f=(r,s)$ with $r>0$.
Then it follows that $\deg (\ad f)^{2m+1}(x)=((2m+1)r+c,(2m+1)(s-1)+d)$.
But each remaining term in the equation for $F_m(\ad f,l_f)(x)$ is of strictly smaller degree.
Hence $P(\ad f,l_f)(a)\neq 0$, which contradicts the assumption on $f$.
\end{proof}

Our approach here is similar to that of \cite{bav-jord} in that we exploit the Poisson structure on $\gr_{lim} A$ to pin down the possible images of the minimal degree element $u\in A$.
However, $u$ is not {\it strictly semisimple} in the sense of \cite[3.3]{bav-jord}.
To determine all possible isomorphisms $D(Q_2,\gamma_2)\rightarrow D(Q_1,\gamma_1)$ we carry out a case-by-case study of the possible images of the standard generators for $D(Q_2,\gamma_2)$.

Hence let $Q_2$ (resp. $Q_1$) be monic of degree $N\geq 3$ (resp. $n\geq 3$) and let $f,g,h$ be the respective images of the standard generators for $D(Q_2,\gamma_2)$ in $D(Q_1,\gamma_1)$.
Assume until further notice that $f=\xi v^r+$ lower terms.
Recall from the proof of Lemma \ref{Fprod} that $F_1(\ad f,l_f)(g)=2\gamma_2$.
Thus $F_1(\ad f,l_f)(2h)=F_1(\ad f,l_f)([f,g])=[f,F_1(\ad f,l_f)(g)]=0$.
It follows that either $g=\xi' wv^{s-1}+$lower terms, or $g=\xi' v^s+$lower terms, for some $\xi'\neq 0$ and $s$.
Similarly, either $h=\xi'' wv^{t-1}+$lower terms, or $h=\xi'' v^t+$lower terms, for some $\xi''\neq 0$ and $t$.
By considering the equalities $[f,g]=2h$, $[f,h]=-2fg+2h+\gamma_2$ and $Q_2(f)+fg^2+h^2=2hg+\gamma_2 g$, we obtain the following exclusive list of possibilities:

\vspace{0.3cm}
(i) $g=\xi' wv^{(N/2-1)r-1}+$lower terms, $h=\xi'' v^{Nr/2}+$lower terms, where $\xi''^2+\xi^N=0$,

(ii) $g=\xi' v^{(N-1)r/2}+$lower terms, $h=\xi'' wv^{(N-1)r/2-1}+$lower terms, where $\xi'^2+\xi^{N-1}=0$,

(iii) $g=\xi' v^{(N-1)r/2}+$lower terms, $h=\xi'' v^t+$lower terms, where $(N-1)r/2<t<Nr/2$ and $\xi'^2+\xi^{N-1}=0$,

(iv) $g=\xi' v^s+$lower terms, $h=\xi'' v^{Nr/2}+$lower terms, where $(N/2-1)r<s<(N-1)r/2$ and $\xi''^2+\xi^N=0$,

(v) $g=\xi' v^s+$lower terms, $h=\xi'' v^{s+r/2}+$lower terms, where $s>(N-1)r/2$ and $\xi\xi'^2+\xi''^2=0$,

(vi) $g=\xi' v^{(N-1)r/2}+$lower terms, $h=\xi'' v^{Nr/2}+$lower terms, where $\xi^N+\xi\xi'^2+\xi''^2=0$.

\vspace{0.3cm}
To deal with these cases, we examine in detail the monomials in $f,g,h$ of highest degree in the expression for $\prod_{i=0}^{m-1} F_i(\ad f,l_f)(g^m)$, and similarly for $hg^{m-1}$.
We will show that the degree of any expression in $f,g,h$ is too high to be equal to $u$ unless $N=3$, where the only possible case is (ii) with $r=1$.

From now on, all monomials in $f,g,h$ will be of the form $f^ih^\epsilon g^j$ with $\epsilon\in\{ 0,1\}$.
For each monomial $x$ in $f,g,h$ and for each non-negative integer $r$, let $J_r(x)$ denote the (finite-dimensional) subspace of $A$ generated by all monomials $f^ih^\epsilon g^j$ with $j+\epsilon\leq r$ and $\deg f^ih^\epsilon g^j<\deg x$.

\begin{lemma}\label{commprep}
Let $f,g,h$ be as in one of the cases (i)-(vi) above, and let $m\geq 2$.

(a) $[f,g^m]=2mhg^{m-1}+m(m-1)g^m+m(m-1)Nf^{N-1}g^{m-2}+a$, for some $a\in J_{m-2}(f^{N-1}g^{m-2})$,

(b) $[f,hg^{m-1}]=-2mfg^m-2(m-1)f^Ng^{m-2}+m(m+1)hg^{m-1}+(2m-1)\gamma' g^{m-1}+(m-1)(m-2)Nf^{N-1}hg^{m-3} +a'$, where $a'\in J_{m-2}(f^{N-1}hg^{m-3})$.
\end{lemma}

\begin{proof}
Clearly $$[f,g^m]= 2\sum_{j=0}^{m-1} g^jhg^{m-1-j} = 2mhg^{m-1}+\sum_{j=0}^{m-1}[g^j,h]g^{m-1-j}$$ and $[g^j,h]=\sum_{l=0}^{j-1}g^l[g,h]g^{j-1-l}$.
Moreover, $g^l[g,h]g^{j-1-l}=g^{j+1}+Ng^lf^{N-1}g^{j-1-l}=g^{j+1}+Nf^{N-1}g^{j-1}+b_l$ for some $b_l\in J_{j-1}(f^{N-1}g^{j-1})$.
It follows that $[g^j,h]=jg^{j+1}+jNf^{N-1}g^{j-1}+b$ for some $b\in J_{j-1}(f^{N-1}g^{j-1})$.
But then clearly $bg^{m-1-j}\in J_{m-2}(f^{N-1}g^{m-2})$.
We deduce that $[f,g^m]\equiv 2mhg^{m-1}+m(m-1)g^m+m(m-1)Nf^{N-1}g^{m-2}\; (\modulo J_{m-2}(f^N g^{m-2}))$.

For (b), $[f,hg^{m-1}]=[f,h]g^{m-1}+h[f,g^{m-1}]$.
By the definition of $D(Q_2,\gamma_2)$, $[f,h]g^{m-1}=-2fg^m+2hg^{m-1}+\gamma_2 g^{m-1}$.
Moreover, $$h[f,g^{m-1}]=2(m-1)h^2g^{m-2}+(m-1)(m-2)hg^{m-1}+(m-1)(m-2)Nf^{N-1}hg^{m-3}+b'$$for some $b'\in hJ_{m-3}(f^{N-1}g^{m-3})\subseteq J_{m-2}(f^{N-1}hg^{m-3})$.
The result now follows from the equality $h^2=-Q_2(f)-fg^2+2hg+\gamma_2 g$.
\end{proof}

\begin{corollary}\label{Fiform}
(a) If $m\geq 2$ then there exists $a\in J_m(fg^m)+J_{m-2}(f^Ng^{m-2})$ such that $$F_i(\ad f,l_f)(g^m)=-4(m^2-i^2)fg^m-4m(m-1)f^Ng^{m-2}+a$$

(b) If $m\geq 3$ then there exists $a'\in J_m(fhg^{m-1})+J_{m-2}(f^Nhg^{m-3})$ such that $$F_i(\ad f,l_f)(hg^{m-1})=-4(m^2-i^2)fhg^{m-1}-4(m-1)(m-2)f^Nhg^{m-3}+a'$$
\end{corollary}

\begin{proof}
By Lemma \ref{commprep}: $$[f,[f,g^m]]=[f,2mhg^{m-1}+m(m-1)g^m+m(m-1)Nf^{N-1}g^{m-1}+a_0]$$ where $a_0\in J_{m-2}(f^Ng^{m-2})$.
But then clearly $[f,g^m]\in J_m(fg^m)$ and $[f,f^{N-1}g^{m-2}],[f,a_0]\in J_{m-2}(f^Ng^{m-2})$.
Applying Lemma \ref{commprep} again, we see that $[f,hg^{m-1}]= -2mfg^m-2(m-1)f^N g^{m-2}+a_1$ for some $a_1\in J_m(fg^m)+J_{m-2}(f^N g^{m-2})$.
Hence the result for $F_i(\ad f,l_f)(g^m)$ follows.

Similarly, Lemma \ref{commprep} implies that \begin{eqnarray}\nonumber\lefteqn{[f,[f,hg^{m-1}]]=[f,-2mfg^m-2(m-1)f^Ng^{m-2}} \\ \nonumber & +m(m+1)hg^{m-1}+(2m+1)\gamma' g^{m-1}+(m-1)(m-2)Nf^{N-1}hg^{m-3}+a_2]\end{eqnarray} where $a_2\in J_{m-2}(f^{N-1}hg^{m-3})$.
But it is immediate that $[f,hg^{m-1}],[f,g^{m-1}]\in J_m(fhg^{m-1})$ and $[f,f^{N-1}hg^{m-3}],[f,a_2]\in J_{m-2}(f^N hg^{m-3})$.
Hence the result for $F_i(\ad f,l_f)(hg^{m-1})$ follows by Lemma \ref{commprep}(a).
\end{proof}

For ease of notation, let $P_i=\prod_{j=0}^i F_j(\ad f,l_f)$ for the rest of this section.
Corollary \ref{Fiform} allows us to describe the monomials in $f,g,h$ which are of highest degree in the expression for $P_i(g^m),P_i(hg^{m-1})$.
We begin with cases (i) and (iv) listed after Lemma \ref{reducetov}.
Here we use the notation $x=\chi f^ih^\epsilon g^j+$lower terms to mean that $x=\chi f^ih^\epsilon g^j+a$, where $a$ is a sum of monomials in $f,g,h$ each of lower degree than $f^ih^\epsilon g^j$.

\begin{lemma}\label{caseiandiv}
Suppose $h=\xi'' v^{Nr/2}+$lower terms and either $g=\xi' wv^{(N/2-1)r-1}+$lower terms (case (i)) or $g=\xi' v^s+$lower terms, where $(N/2-1)r<s<(N-1)r/2$ (case (iv)).
Then for any $m\geq 1$:

(a) $P_i(g^{2m})=\left\{\begin{array}{ll} \chi_{i} f^{iN}hg^{2m-2i-1}+\mbox{lower terms} & \mbox{if $0\leq i<m$,} \\
\chi_{i} f^{i+(m-1)(N-1)}hg+\mbox{lower terms} & \mbox{if $m\leq i<2m$.} \end{array}\right.$

(b) $P_i(g^{2m-1})=\left\{\begin{array}{ll} \chi_{i} f^{iN}hg^{2m-2i-2}+\mbox{lower terms} & \mbox{if $0\leq i<m$,} \\
\chi_{i} f^{i+(m-1)(N-1)}h+\mbox{lower terms} & \mbox{if $m\leq i<2m-1$.}\end{array}\right.$

(c) $P_i(hg^{2m-1})=\left\{\begin{array}{ll} \eta_{i} f^{(i+1)N}g^{2m-(2i+2)}+\mbox{lower terms} & \mbox{if $0\leq i<m$,} \\
\eta_{i} f^{i+1+m(N-1)}+\mbox{lower terms} & \mbox{if $m\leq i<2m$.}\end{array}\right.$

(d) $P_i(hg^{2m-2})=\left\{\begin{array}{ll} \eta_{i} f^{(i+1)N}g^{2m-(2i+3)}+\mbox{lower terms} & \mbox{if $0\leq i<m-1$,} \\
\eta_{i} f^{i+1+(m-1)(N-1)}g+\mbox{lower terms} & \mbox{if $m-1\leq i<2m-1$.}\end{array}\right.$

Here $\chi_{i}$ (resp. $\eta_{i}$) is a real number of sign $(-1)^i$ (resp. $(-1)^{i+1}$).
\end{lemma}

\begin{proof}
Our proof is by induction on $m$ and $i$.
Since $P_0=\ad f$, the lemma is true for $i=0$ by Lemma \ref{commprep} and the fact that $\deg f^N>\deg fg^2$.
For $m=1$, this proves (b) and (d).
By a direct calculation, $P_1(g^2)=-48fhg+$lower terms and $P_1(hg)=24f^{N+1}+$lower terms.
Hence (a) and (c) are also true for $m=1$.
We assume therefore that $m\geq 2$.

By Cor. \ref{Fiform}, $P_1(g^l)=[f,-4(l^2-1)fg^l-4l(l-1)f^Ng^{l-2}+a]$ for some $a\in J_l(fg^l)+J_{l-2}(f^Ng^{l-2})$.
Let $\delta$ be equal to $(r,-1)$ in case (i), and equal to $(Nr/2-s,0)$ in case (iv).
Then $\deg [f,g^l]=\deg g^l+\delta$ for any $l\geq 1$.
Moreover, if $x$ is any monomial in $f,g,h$ and $\sum a_{ij\epsilon} f^ih^\epsilon g^j$ is the unique expression for $[f,x]$ in terms of monomials in $f,g,h$ then it follows from Lemma \ref{commprep} that each non-zero term $a_{ij\epsilon}f^ih^\epsilon g^j$ has degree less than or equal to $\deg x+\delta$.
But therefore $[f,a]$ and $[f,fg^l]$ are both of degree less than $f^Nhg^{l-3}$.
Hence $P_1(g^l)=-8l(l-1)(l-2)f^Nhg^{l-3}+$lower terms for any $l\geq 3$.
This proves (a) and (b) for $i=1$.
A direct calculation establishes that $P_1(hg^2)=144f^{N+1}g+$lower terms.
Hence (d) is true for $m=2$ and $i=1$.
We therefore consider $P_1(hg^{l-1})$ for $l\geq 4$.
By Cor. \ref{Fiform}, $F_1(\ad f,l_f)(hg^{l-1})=-4(l^2-i^2)fhg^{l-1}-4(l-1)(l-2)f^Nhg^{l-3}+a'$, where $a'\in J_l(fhg^{l-1})+J_{l-2}(f^Nhg^{l-3})$.
The highest degree term here is $f^Nhg^{l-3}$.
Moreover, $[f,f^Nhg^{l-3}]=-2(l-3)f^{2N}g^{l-4}+$lower terms, and $\deg f^{2N}g^{l-4}=\deg f^Nhg^{l-3}+\delta$.
By the remarks above, $P_1(hg^{l-1})=8(l-1)(l-2)(l-3)f^{2N}g^{l-4}+$lower terms, which confirms (c) and (d) for $i=1$.

An equivalent statement for the Lemma can be formulated in terms of degrees (and leading coefficients) of the $P_i(g^l),P_i(hg^{l-1})$.
Specifically: $$\deg P_i(g^{2m})-\deg g^{2m}=\left\{\begin{array}{ll}
(2i+1)\delta & \mbox{if $i<m$,} \\
(2m-1)\delta+(i-m+1)(r,0) & \mbox{if $m\leq i<2m$.}\end{array}\right.$$
$$\deg P_i(g^{2m-1})-\deg g^{2m-1}=\left\{\begin{array}{ll} (2i+1)\delta & \mbox{if $i<m$,} \\
(2m-1)\delta+(i-m+1)(r,0) & \mbox{if $m\leq i<2m-1$.}\end{array}\right.$$
$$\deg P_i(hg^{2m-1})-\deg hg^{2m-1}=\left\{\begin{array}{ll}
(2i+1)\delta & \mbox{if $i<m$,} \\
(2m-1)\delta+(i-m+1)(r,0) & \mbox{if $m\leq i<2m$.}\end{array}\right.$$
$$\deg P_i(hg^{2m-2})-\deg hg^{2m-2}=\left\{\begin{array}{ll}
(2i+1)\delta & \mbox{if $i<m-1$,} \\
(2m-3)\delta+(i-m+2)(r,0) & \mbox{if $m-1\leq i<2m-1$.}\end{array}\right.$$

(We retain of course the assumption on the signs of the leading coefficients $\chi_i,\eta_i$.)

Assume therefore that $i\geq 2$ and that (a)-(d) are known to be true for all pairs $(m',i')$ with $m'<m$ or $m'=m$ and $i'<i$.
By Cor. \ref{Fiform}, $F_i(\ad f,l_f)(g^{2m})=-4(4m^2-i^2)fg^{2m}-8m(2m-1)f^Ng^{2m-2}+a$ for some $a\in J_{2m}(fg^{2m})+J_{2m-2}(f^Ng^{2m-2})$.
But let $a=a_1+a_2$, where $a_1\in J_{2m}(fg^{2m})$ and $a_2\in J_{2m-2}(f^Ng^{2m-2})$.
By the induction hypothesis and the remarks above, $\deg P_{i-1}(a_1)-\deg a_1\leq \deg P_{i-1}(fg^{2m})-\deg fg^{2m}$ and $\deg P_{i-1}(a_2)-\deg a_2\leq \deg P_{i-1}(f^Ng^{2m-2})-\deg f^Ng^{2m-2}$.
Hence $P_i(g^{2m})=-4(4m^2-i^2)P_{i-1}(fg^{2m})-8m(2m-1)P_{i-1}(f^Ng^{2m-2})+$lower terms.
If $i<m$, then by the induction hypothesis $P_{i-1}(fg^{2m})=\chi_{i-1}f^{N(i-1)+1}hg^{2m-2i+2}+$lower terms and $P_{i-1}(f^Ng^{2m-2})=\chi_{i-1}'f^{iN}hg^{2m-2i-1}+$lower terms, where $\chi_{i-1}$ and $\chi'_{i-1}$ are both of sign $(-1)^{i-1}$.
It follows that $P_i(g^{2m})=-8m(2m-1)\chi_{i-1}f^{iN}hg^{2m-2i-1}+$lower terms, which proves the induction step for (a) in the case $i<m$.
If $2m-1>i\geq m$, then by the induction hypothesis $P_{i-1}(fg^{2m})=\chi_{i-1}f^{i+(m-1)(N-1)}hg+$lower terms and $P_{i-1}(f^Ng^{2m-2})=\chi_{i-1}'f^{i+(m-1)(N-1)}hg+$lower terms, where $\chi_{i-1}$ and $\chi'_{i-1}$ are of sign $(-1)^{i-1}$.
This proves the induction step in this case.
Finally, $P_{2m-2}(f^Ng^{2m-2})=P_{2m-2}(a_2)=0$, hence $P_{2m-1}(g^{2m})=-4(4m-1)P_{2m-2}(fg^{2m})+$lower terms.
But $P_{2m-2}(g^{2m})=\chi_{2m-2}f^{2m-1+(m-1)(N-1)}hg+$lower terms, where $\chi_{2m-2}$ is positive.
It follows that $P_{2m-1}(g^{2m})=-4(4m-1)\chi_{2m-2}f^{m+(m-1)N}hg+$lower terms.

This proves the induction step for (a).
The arguments for (b) and (c) are identical.
For (d) we need to be slightly careful, for if $x$ is a monomial in $J_{2m-1}(fhg^{2m-2})$ then it is not necessarily true that $\deg P_{i-1}(x)-\deg x\leq \deg P_{i-1}(fhg^{2m-2})-\deg fhg^{2m-2}$ (and similarly for $J_{2m-3}(f^Nhg^{2m-4})$).
In fact, one can see easily from the description of degrees above that if $x$ is a monomial in $J_{2m-1}(fhg^{2m-2})$ then $\deg P_{i-1}(x)-\deg x>\deg P_{i-1}(fhg^{2m-2})-\deg fhg^{2m-2}$ if and only if $i=m$ and $x=g^{2m-1}$ or $x=fg^{2m-1}$.
However, in this case we still have that $\deg P_{i-1}(x)\leq \deg P_{i-1}(fhg^{2m-2})$.
Similarly, if $x\in J_{2m-3}(f^Nhg^{2m-4})$ then $\deg P_{i-1}(x)\leq \deg P_{i-1}(f^Nhg^{2m-4})$.
It follows that $P_i(hg^{2m-2})=-4(4m^2-i^2)P_{i-1}(fhg^{2m-2})-4(m-1)(m-2)P_{i-1}(f^Nhg^{2m-4})+$lower terms.
The rest of the argument now proceeds as above.
\end{proof}

\begin{corollary}\label{iiv}
Suppose $g,h$ are as in Lemma \ref{caseiandiv}.
Then there is no possible expression for $u$ in terms of $f,g,h$.
\end{corollary}

\begin{proof}
Suppose there exists such an expression $u=\sum_{i,j\geq 0,\epsilon\in\{ 0,1\} } a_{ij\epsilon} f^ih^\epsilon g^j$, and let $m=\max_{\{ a_{ij\epsilon}\neq 0\}}(j+\epsilon)$.
Clearly $\deg g^m<\deg hg^{m-1}<\deg fg^m$.
Moreover, $\deg P_{m-1}(g^m)<\deg P_{m-1}(hg^{m-1})<\deg P_{m-1}(fg^m)$ by Lemma \ref{caseiandiv}.
Applying $P_{m-1}$ to both sides of the equation, we have the equality $P_{m-1}(u)=\sum_{j+\epsilon=m}a_{ij\epsilon}P_{m-1}(f^ih^\epsilon g^j)$.
Thus $\deg P_{m-1}(u)\geq \deg P_{m-1}(g^m)$.
To show that there can be no such expression for $u$, it will therefore suffice to show that $\deg P_{m-1}(u)<\deg P_{m-1}(g^m)$.
If $m=1$, then $\deg P_{m-1}(u)=(r,1)<\deg P_{m-1}(g)=(Nr/2,0)$.
Suppose therefore that $m\geq 2$, hence $\deg P_{m-1}(u)<((2m-1)r,0)$.

If $m$ is even, then by Lemma \ref{caseiandiv}, $P_{m-1}(g^m)=\chi_{m-1}f^{(m/2-1)(N+1)+1}hg+$lower terms, hence $\deg P_{m-1}(g^m)=(m(N+1)r/2,1)$ in case (i) and $\deg P_{m-1}(g^m)=(m(N+1)r/2+(s-(N/2-1)r),0)>(m(N+1)r/2,0)$ in case (iv).
But $N+1\geq 4$, hence $\deg P_{m-1}(g^m)>((2m-1)r,0)$ in both cases.
Similarly, if $m$ is odd then $P_{m-1}(g^m)=\chi_{m-1}f^{(m-1)(N+1)/2}h+$lower terms, hence $\deg P_{m-1}(g^m)=((m(N+1)-1)r/2,0)>((2m-1)r,0)$.
This completes the proof.
\end{proof}

Next we deal with case (v).
This case is fairly straightforward, since the highest degree term in the expression for $h^2$ is $fg^2$.
Once more we write $x=\chi f^ih^\epsilon g^j+$lower terms to mean $x=\chi f^ih^\epsilon g^j+a$, where $a$ is a sum of terms $a_{kl\epsilon}f^kh^\epsilon g^l$, each of degree strictly less than that of $f^ih^\epsilon g^j$.

\begin{lemma}\label{casev}
Suppose $g=\xi' v^s+$lower terms and $h=\xi'' v^{s+r/2}+$lower terms, where $s>(N-1)r/2$ (hence $r$ is even).
Then for any $0\leq i<m$, $P_i(g^m)=\chi_{i}f^ihg^{m-1}+$lower terms, where $\chi_{i}$ is a real number of sign $(-1)^i$ and $P_i(hg^{m-1})=\eta_{i}f^{i+1}g^m+$lower terms, where $\eta_{i}$ is a real number of sign $(-1)^{i+1}$.
\end{lemma}

\begin{proof}
We apply a similar argument to that in the proof of Lemma \ref{caseiandiv}.
The statement of the Lemma for $i=0$ follows immediately from Lemma \ref{commprep}.
Hence assume $i\geq 1$, and that the Lemma is known to be true for all pairs $(m',i')$ with $i'<m'$ and either $m'<m$ or $m'=m$ and $i'<i$.
Note that the induction hypothesis implies that $\deg P_{i'}(g^{m'})-\deg g^{m'}=(2i'+1)(r/2,0)$ for any such pair $(m',i')$, and similarly for $hg^{m-1}$.
It follows that if $a\in J_m(fg^m)$ (resp. $a'\in J_m(fhg^{m-1})$) then $\deg P_{i-1}(a)\leq \deg f^ihg^{m-1}$ (resp. $\deg P_{i-1}(a')\leq\deg f^{i+1}g^m$).
By Lemma \ref{Fiform}, $F_i(\ad f,l_f)(g^m)=-4(m^2-i^2)fg^m+a$, where $a\in J_m(fg^m)$.
Moreover, $P_{i-1}(fg^m)=\chi_{i-1}f^ihg^{m-1}+$lower terms, where $\chi_{i-1}$ is a real number of sign $(-1)^{i-1}$.
Hence $P_i(g^m)=P_{i-1}(-4(m^2-i^2)fg^m+a)=-4(m^2-i^2)\chi_{i-1}f^ihg^{m-1}+$lower terms.
This proves the induction step for $P_i(g^m)$.
Similarly, $F_i(\ad f,l_f)(hg^{m-1})=-4(m^2-i^2)fhg^{m-1}+a'$, where $a'\in J_m(fhg^{m-1})$.
But $P_{i-1}(fhg^{m-1})=\eta_{i-1}f^{i+1}g^m+$lower terms, where $\eta_{i-1}$ is real of sign $(-1)^i$.
Since $P_{i-1}(a')\in J_m(f^{i+1}g^m)$, we deduce that $P_i(a')=-4(m^2-i^2)f^{i+1}g^m+$lower terms.
This completes the proof.
\end{proof}

\begin{corollary}\label{v}
Suppose $g$ and $h$ are as in Lemma \ref{casev}.
Then there is no possible expression for $u$ in terms of $f,g,h$.
\end{corollary}

\begin{proof}
Suppose there exists an expression $u=\sum a_{ij\epsilon} f^ih^\epsilon g^j$, and as in the proof of Lemma \ref{caseiandiv}, let $m=\max_{\{ a_{ij\epsilon}\neq 0\} }(j+\epsilon)$.
Applying $P_{m-1}$ to both sides, we have an equality $P_{m-1}(u)=\sum_{j+\epsilon=m}a_{ij\epsilon}f^i P_{m-1}(h^\epsilon g^j)$.
Moreover, it is immediate from Lemma \ref{casev} that $\deg P_{m-1}(g^m)<\deg P_{m-1}(hg^{m-1})<\deg fP_{m-1}(g^m)$, hence $\deg P_{m-1}(u)\geq \deg P_{m-1}(g^m)$.
Hence to prove the lemma we have only to prove that $\deg P_{m-1}(u)<\deg P_{m-1}(g^m)$.
Clearly $[u,v^r]=2r wv^{r-1}+$lower terms, hence $\deg P_0(u)=(r,1)$.
But $P_0(g)=2h$ is of degree $(s+r/2,0)> (Nr/2,0)>(r,1)$.
On the other hand, if $m\geq 2$ then $\deg P_{m-1}(u)< ((2m-1)r,0)$.
Furthermore, $\deg P_{m-1}(g^m)=(ms+(m+1/2)r,0)>((m(N+1)-1)r/2,0)$.
Since $N\geq 3$, $\deg P_{m-1}(g^m)> \deg P_{m-1}(u)$.
This completes the proof.
\end{proof}

We have therefore eliminated cases (i), (iv) and (v) listed after Lemma \ref{reducetov}.
Roughly speaking, the highest degree monomial in the expression for $h^2$ ($f^N$ in cases (i) and (iv), $fg^2$ in case (v)) contributes the highest degree monomial in the expression for $F_i(\ad f,l_f)(g^m)$, and thus eventually in the expression for $P_i(g^m)$ (and similarly $P_i(hg^{m-1})$).
For the remaining cases, we must replace $h^2$ by terms of possibly higher degree, since we wish to find the expressions for $P_{m-1}(g^m)$ and $P_{m-1}(hg^{m-1})$ in terms of the monomials $f^ih^\epsilon g^j$ with $\epsilon\in\{ 0,1\}$.
In these circumstances $f^N$ and $fg^2$ are now of equal degree, hence our final expression for the leading term of $P_{m-1}(g^m)$ will contain a number of monomials in $f,g,h$ of equal degree.
Here we write $x=\sum \chi_j f^{i+j(N-1)}hg^{m-2j-1}+$lower terms (resp. $x=\sum \eta_j f^{i+1+j(N-1)}g^{m-2j}+$lower terms) to mean that $x=\sum \chi_j f^{i+j(N-1)}hg^{m-2j-1}+a$ (resp. $x=\sum \eta_j f^{i+1+j(N-1)}g^{m-2j}+a$), where $a$ is a sum of monomials in $f,g,h$, each of degree less than that of $f^ihg^{m-1}$ (resp. $f^{i+1}g^m$).

\begin{lemma}\label{othercases}
Suppose $g=\xi' v^{(N-1)r/2}+$lower terms and either $h=\xi'' wv^{(N-1)r/2-1}+$lower terms (case (ii)) or $h=\xi' v^t+$lower terms, where $(N-1)r/2<t\leq Nr/2$ (cases (iii) and (vi)).
Then for any $i<m$:

(a) Let $l=\min\{i,[(m-1)/2]\}$.
Then $P_i(g^m)=\sum_{j=0}^{l} \chi_{j} f^{i+j(N-1)}hg^{m-2j-1}+$lower terms, where the $\chi_{j}$ are real numbers of sign $(-1)^i$.

(b) Let $l=\min\{i+1,[m/2]\}$.
Then $P_i(hg^{m-1})=\sum_{j=0}^{l} \eta_{j} f^{i+1+j(N-1)}g^{m-2j}+$lower terms, where the $\eta_{j}$ are real numbers of sign $(-1)^{i+1}$.
\end{lemma}

\begin{proof}
We follow a similar argument to the proofs of Lemmas \ref{caseiandiv} and \ref{casev}.
If $i=0$, then (a) and (b) are direct consequences of Lemma \ref{commprep}.
Assume therefore that $m>i\geq 1$ and that the Lemma is known to be true for all pairs $(m',i')$ with $i'<m'$ and either $m'<m$ or $m'=m$ and $i'<i$.
By a direct calculation, $P_1(\ad f,l_f)(g^2)=-48fhg+$lower terms and $P_1(hg)=48f^2g^2+24f^{N+1}+$lower terms.
Hence assume $m\geq 3$.
Let $\delta=\deg h-\deg g\leq (r/2,0)$.
By Cor. \ref{Fiform}, $F_i(\ad f,l_f)(g^m)=-4(m^2-i^2)fg^m-4m(m-1)f^Ng^{m-2}+a$ and $F_i(\ad f,l_f)(hg^{m-1})=-4(m^2-i^2)fhg^{m-1}-4(m-1)(m-2)f^Nhg^{m-3}+a'$, where $a\in J_m(fg^m)$ and $a'\in J_m(f^Nhg^{m-1})$.
We note that if $i-1<m'\leq m$ then by the induction hypothesis $P_{i-1}(g^{m'})$ is a sum of monomials in $f,g,h$, each of degree less than or equal to $\delta+(ir,0)+\deg g^m$.
On the other hand, $P_{i-1}(hg^{m'-1})$ is a sum of monomials of degree less than or equal to $((i+1)r,0)-\delta+\deg hg^{m'-1}$.
But therefore $P_{i-1}(a')\in J_m(f^{i+1}g^m)$ for any $a'\in J_m(hg^{m-1})$.
It follows by the induction hypothesis that $$\begin{array}{rl}P_i(hg^{m-1})=& P_{i-1}(-4(m^2-i^2)fhg^{m-1}-4(m-1)(m-2)f^Nhg^{m-3}+a') \\ = & -4(m^2-i^2)P_{i-1}(fhg^{m-1})-4(m-1)(m-2)P_{m-1}(f^Nhg^{m-3})+b'\end{array}$$
where $b'\in J_m(f^{i+1}g^m)$.
This proves the induction step for (b).

For (a) we have to be careful, for it is not in general true that if $x$ is a monomial in $f,g,h$ then $\deg P_{i-1}(x)-\deg x\leq \deg P_{i-1}(g^m)-\deg g^m$.
In fact we can see that this is true if and only if $x=f^{i'}g^{m'}$ for some $i',m'$.
On the other hand, if $x\in J_m(fg^m)$ is of the form $f^{i'}hg^{m'-1}$ then $\deg x\leq\deg g^m-(r,0)+\delta$.
By our statement above on the degrees of $P_{i-1}(g^{m'})$, $P_{i-1}(hg^{m'-1})$ we nevertheless have that $P_{m-1}(x)\in J_m(f^ihg^{m-1})$, hence that $P_{m-1}(a)\in J_m(f^ihg^{m-1})$.
The argument now proceed exactly as above.
\end{proof}

Suppose therefore that $g,h$ are as in Lemma \ref{othercases}.
Then $P_{m-1}(g^m)=\chi_{0}f^{m-1}hg^{m-1}+\chi_{1}f^{m+N-2}hg^{m-3}+\ldots+a$, where $a\in J_m(f^{m-1}hg^{m-1})$ and the $\chi_{j}$ are real numbers of the same sign.
Since $\deg g^2=\deg f^{N-1}$, the monomials $f^{m-1+j(N-1)}hg^{m-2j-1}$ are of equal degree.
We ask therefore whether it is possible that the highest degree terms of these monomials (expressed in terms of $u,v$, and $w$) cancel out.
Specifically, this holds if and only if $\chi_{0}\xi'^{(m-1)}+\chi_{1}\xi^{N-1}\xi'^{(m-3)}+\chi_{2}\xi^{2(N-1)}\xi'^{(m-5)}+\ldots=0$, which is in turn true if and only if $\chi_{0}+\chi_{1}\mu+\chi_{2}\mu^2+\ldots=0$, where $\mu=\xi^{N-1}/(\xi')^{2}$.
In case (ii) or (iii), $\mu=-1$.
In case (vi) $\mu\neq -1$, since $\xi^N+\xi\xi'^2+\xi''^2=0$.

\begin{lemma}\label{coeffs}
(a) $P_{m-1}(g^m)=\chi_0 f^{m-1}hg^{m-1}+\chi_1 f^{N+m-2}hg^{m-3}+\ldots+a$, where $a\in J_m(f^{m-1}hg^{m-1})$ and $\chi_i=\chi_0\cdot\left(\begin{array}{c} m-i-1 \\ i\end{array}\right)/4^i$ for $0\leq i\leq [(m-1)/2]$.

(b) $P_{m-1}(hg^{m-1})=\eta_0 f^mg^m+\eta_1 f^{m+N-1}g^{m-2}+\ldots+a'$, where $a'\in J_m(f^mg^m)$ and $\eta_i=\eta_0\cdot(\left(\begin{array}{c} m-i \\ i\end{array}\right)+\left(\begin{array}{c} m-i-1\\ i-1 \end{array}\right) )/4^i$ for $0\leq i\leq [m/2]$.
\end{lemma}

\begin{proof}
The fact that $P_{m-1}(g^m)$ has the above form for some constants $\chi_0,\chi_1,\ldots$ follows immediately from Lemma \ref{othercases}.
Moreover, $F_m(\ad f,l_f)(P_{m-1}(g^m))=0$.
Let $\omega_i=(m-2i+1)(m-2i)\chi_{i-1}+((m-2i)^2-m^2)\chi_i$ for $1\leq i\leq [(m-1)/2]$.
By application of (the argument in the proof of) Cor. \ref{Fiform}, $$F_m(\ad f,l_f)(P_{m-1}(g^m))=-4\sum_1^{[(m-1)/2]}\omega_if^{m-1+i(N-1)}hg^{m-2i-1}+F_m(\ad f,l_f)(a)$$
But by the observation in the proof of Lemma \ref{othercases}, $F_{m}(\ad f,l_f)(a)\in J_m(f^{N+m-1}hg^{m-3})$.
Hence each of the coefficients $(m-2i)(m-2i+1)\chi_{i-1}-4i(m-i)\chi_i$ is equal to zero.
We deduce that $$\frac{\chi_i}{\chi_0}=\frac{(m-1)!}{(m-2i-1)!}\cdot\frac{1}{i!}\cdot\frac{(m-i-1)!}{(m-1)!}\cdot \frac{1}{4^i}=\left(\begin{array}{c} m-i-1 \\ i\end{array}\right)/4^i$$

Similarly, the existence of some constants $\eta_i$ and an expression for $P_{m-1}(hg^{m-1})$ as in (b) follows immediately from Lemma \ref{othercases}.
We apply the same argument as above.
Thus let $\omega'_i= (m-2i+2)(m-2i+1)\eta_{i-1}-4i(m-i)\eta_i$ for $1\leq i\leq [m/2]$.
Then $$F_m(\ad f,l_f)(P_{m-1}(hg^{m-1}))=-4\sum_0^{[m/2]} \omega_i'f^{m+1+i(N-1)}g^{m-2i}+F_m(\ad f,l_f)(a')=0$$
By Lemma \ref{othercases}, $F_m(\ad f,l_f)(a')\in J_m(f^{m+N}g^{m-2})$.
Hence each of the coefficients $(m-2i+2)(m-2i+1)\eta_{i-1}-4i(m-i)\eta_i$ is equal to zero.
We conclude that $$\frac{\eta_i}{\eta_0}=\frac{m!}{(m-2i)!}\cdot\frac{(m-i-1)!}{i!(m-1)!}\cdot\frac{1}{4^i}=\frac{(m-i-1)!}{i!(m-2i)!}\cdot\frac{m}{4^i}$$
from which (b) follows.
\end{proof}

We thus introduce the polynomials $$p_m(t)=\sum_0^{[(m-1)/2]}\left(\begin{matrix} m-i-1 \\ i \end{matrix}\right)(t/4)^i$$ and $$q_m(t)=\sum_0^{[m/2]}(\left(\begin{matrix} m-i \\ i \end{matrix}\right)+\left(\begin{matrix} m-i-1 \\ i-1 \end{matrix}\right) )(t/4)^i.$$

\begin{lemma}\label{binomials}
(a) $p_{m+1}(t)=p_m(t)+tp_{m-1}(t)/4$ and $q_{m+1}(t)=q_m(t)+tq_{m-1}(t)$.

(b) $q_m(t)=p_{m}(t)+tp_{m-1}(t)/2$.

(c) $p_m(-1)=m/2^{m-1}$ and $q_m(-1)=1/2^{m-1}$.
\end{lemma}
\begin{proof}
Part (a) follows immediately from the fact that $\left(\begin{matrix} m-i-1 \\ i \end{matrix}\right)-\left(\begin{matrix} m-i-2 \\ iÊ\end{matrix}\right)=\left(\begin{matrix} m-i-2 \\ i-1 \end{matrix}\right)$ and similarly for $\left(\begin{matrix} m-i \\ i \end{matrix}\right)$, $\left(\begin{matrix} m-i-1 \\ i-1 \end{matrix}\right)$.
For part (b), we see by a simple re-indexing exercise that $q_m(t)=p_{m+1}(t)+tp_{m-1}(t)/4$, hence that the equality is true by application of (a).
If $\alpha_m=p_m(-1)$ and $\beta_m=q_m(-1)$ then it follows that $\alpha_{m+1}=\alpha_m-\alpha_{m-1}/4$ and $\beta_{m+1}=\beta_m-\beta_{m-1}/4$.
The general solution to this difference equation is $(Am+B)/2^m$.
But $\alpha_2=1,\alpha_3=3/4$.
Hence $\alpha_m=m/2^m$.
Similarly, $\beta_2=1/2$ and $\beta_3=1/4$, hence $\beta_m=1/2^{m-1}$.
\end{proof}

\begin{corollary}\label{notzero}
If $h=\xi'' wv^{(N-1)r/2-1}+$lower terms (case (ii)) or $h=\xi'' v^t+$lower terms, where $(N-1)r/2<t<Nr/2$ (case (iii)) then $\deg P_{m-1}(g^m)=\deg f^{m-1}hg^{m-1}$ and $\deg P_{m-1}(hg^{m-1})=\deg f^mg^m$.
\end{corollary}

\begin{proof}
We remarked after Lemma \ref{othercases} that the degree of $P_{m-1}(g^m)$ (resp. $P_{m-1}(hg^{m-1})$) is lower than that of $f^{m-1}hg^{m-1}$ (resp. $f^mg^m$) if and only if $\chi_0+\chi_1\mu+\ldots+\chi_{[(m-1)/2]}\mu^{[(m-1)/2]}=0$ (resp. $\eta_0+\eta_1\mu+\ldots+\eta_{[m/2]}\mu^{[m/2]}=0$).
But here $\mu=-1$, hence by Lemmas \ref{coeffs} and \ref{binomials} $\chi_0+\chi_1\mu+\ldots+\chi_{[(m-1)/2]}\mu^{[(m-1)/2]}\neq 0$ and $\eta_0+\eta_1\mu+\ldots+\eta_{[m/2]}\mu^{[m/2]}\neq 0$.
\end{proof}

\begin{corollary}\label{notzero2}
If $\mu\neq -1$ then there exists no $m$ such that $p_m(\mu)$ and $q_m(\mu)$ are both zero.
Hence if $h=\xi'' v^{Nr/2}+$lower terms (case (vi)) then either $\deg P_{m-1}(g^m)=\deg f^{m-1}hg^{m-1}$ or $\deg P_{m-1}(hg^{m-1})=\deg f^mg^m$.
\end{corollary}

\begin{proof}
Suppose there exists $m$ such that $p_m(\mu)=p_{m+1}(\mu)=0$.
Let $m'$ be minimal such.
Then by Lemma \ref{binomials}(a) $p_{m'-1}(\mu)=0$, which contradicts the minimality of $m'$.
Hence there exists no $m$ such that $p_m(\mu)=p_{m+1}(\mu)=0$.
But if $p_m(\mu)=q_m(\mu)=0$ then $p_{m-1}(\mu)=0$ by Lemma \ref{binomials}(b).
\end{proof}

\begin{lemma}\label{potential}
Let $a\in J_{m+1}(f^ihg^m)$ and $a'\in J_m(f^ig^m)$.
If $h=\xi'' v^{Nr/2}+$lower terms (case (vi)) then $[f,a]\in J_{m+1}(f^{i+1}g^{m+1})$ and $[f,a']\in J_m(f^ihg^{m-1})$.
\end{lemma}

\begin{proof}
In case (vi), we have $[f,g]=2h$ and $\deg h-\deg g=(r/2,0)$.
Moreover, $[f,h]=-2fg+$lower terms, and $\deg fg-\deg h=(r/2,0)$.
It follows that $\deg [f,y]-\deg y\leq (r/2,0)$ for any monomial $y$ in $f,g,h$.
Hence $[f,a]\in J_{m+1}(f^{i+1}g^{m+1})$ for any $a\in J_{m+1}(f^ihg^m)$ and $[f,a']\in J_m(f^ihg^{m-1})$ for any $a'\in J_m(f^ig^m)$.
This completes the proof.
\end{proof}

Note that Lemma \ref{potential} is not true in cases (ii) and (iii).
(Hence the proof of Lemma \ref{Ngeq4} really requires different arguments for the cases $\mu=-1$, $\mu\neq -1$.)

\begin{lemma}\label{Ngeq4}
Suppose $f,g,h$ are as in Lemma \ref{othercases}.

(a) If $N\geq 4$ then there is no possible expression for $u$ in terms of $f,g,h$.

(b) If $N=3$ then there exists no possible expression for $u$ in terms of $f,g,h$ unless $f=\xi v+$lower terms, $g=\xi' v+$lower terms, $h=\xi'' w+$lower terms.
Moreover, in this case any expression for $u$ must be of the form $c_1 g+c_2 f+c_3$, where $c_1,c_2\in {\mathbb C}^\times$ and $c_3\in{\mathbb C}$.
\end{lemma}

\begin{proof}
Suppose there is such an expression $u=\sum_{ij\epsilon}a_{ij\epsilon}f^ih^\epsilon g^j$ (with the sum taken over all $i,j\geq 0$, $\epsilon\in\{ 0,1\}$).
Let $m=\max_{a_{ij\epsilon}\neq 0}(j+\epsilon)$.
Assume first of all that $m>1$: we will show that such an expression is impossible for all $N$.
Suppose first of all that $a_{im0}f^ig^m$ is the term of highest degree in the expression for $u$ among those of the form $f^jg^m$, $f^jhg^{m-1}$.
By Lemma \ref{othercases}, $P_{m-1}(u)=a_{im0}f^iP_{m-1}(g^m)+a$, where $a\in J_m(f^{i+m-1}hg^{m-1})$.
If $\mu=-1$ (cases (ii) and (iii)) then by Cor. \ref{notzero}, $\deg P_{m-1}(g^m)=\deg f^{m-1}hg^{m-1}$.
But $\deg f^{m-1}hg^{m-1}\geq ((m(N+1)/2-1)r,0)\geq ((2m-1)r,0)>\deg P_{m-1}(u)$, hence such an equality is impossible.
On the other hand, if $\mu\neq -1$ (case (vi)) then either $\deg P_{m-1}(g^m)=\deg f^{m-1}hg^{m-1}$ or $\deg P_{m-1}(hg^{m-1})=\deg f^mg^m$.
If $\deg P_{m-1}(g^m)=\deg f^{m-1}hg^{m-1}$, then the argument above provides a contradiction.
If not, then we consider the equality $[f,u]=2ma_{im0}f^ihg^{m-1}+a'$, where $a'\in J_m(f^ihg^{m-1})$.
By Lemmas \ref{potential} and \ref{othercases}, $P_{m-1}([f,u])=2ma_{im0}f^iP_{m-1}(hg^{m-1})+a$, where $a\in J_m(f^mg^m)$.
But now $\deg P_{m-1}([f,u])<(2mr,0)$ and $\deg f^iP_{m-1}(hg^{m-1})\geq \deg f^mg^m=(m(N+1)r/2,0)$.
Hence there is no such expression for $u$.

Assume therefore that $a_{i(m-1)1}f^ihg^{m-1}$ is the highest degree term in the expression for $u$ among those of the form $f^jg^m$, $f^jhg^{m-1}$.
Once more we apply $P_{m-1}$: $P_{m-1}(u)=a_{i(m-1)1}f^iP_{m-1}(hg^{m-1})+a$, where $a\in J_m(f^{i+m}g^m)$.
If $\mu=-1$ (cases (ii) and (iii)) then $\deg P_{m-1}(u)<((2m-1)r,0)$ and $\deg f^iP_{m-1}(hg^{m-1})\geq (m(N+1)r/2,0)$, hence such an equality is impossible.
If $\mu\neq -1$, then either $\deg P_{m-1}(hg^{m-1})=(m(N+1)r/2,0)$ or $\deg P_{m-1}([f,hg^{m-1}])=((m(N+1)+1)r/2,0)$.
In the first case, the argument for $\mu=-1$ shows that equality of degrees is impossible.
In the second case, by Lemma \ref{potential} there is an equality $[f,u]=-2ma_{i(m-1)1}f^{i+1}g^m+a'$, where $a'\in J_m(f^{i+1}g^m)$.
Then $\deg P_{m-1}([f,u])<(2mr,0)$.
But $\deg f^iP_{m-1}(f^{i+1}g^m)>(m(N+1)r/2,0)$.
This proves our claim.

Suppose therefore that $m=1$.
Hence $u=m_1(f)g+m_2(f)h+m_3(f)$ for some polynomials $m_1(t),m_2(t),m_3(t)$.
Applying $(\ad f)$, we have an equality $[f,u]=m_1(f)[f,g]+m_2(f)[f,h]$.
But $\deg h=\deg [f,g]<\deg [f,h]<\deg f[f,g]$ and $\deg [f,u]=(r,1)$, hence such an equality is only possible if $h=\xi'' w+$lower terms, that is to say, if $N=3$, $r=1$ and we are in case (ii).
Thus (b) follows.
\end{proof}

Lemma \ref{Ngeq4} essentially completes our determination of the isomorphisms $D(Q_2,\gamma_2)\rightarrow D(Q_1,\gamma_1)$ in the case $N\geq 4$ (or $n\geq 4$).
The following straightforward lemma is the final step.

\begin{lemma}\label{outer}
Let $Q_1$ and $Q_2$ be monic polynomials of degree $\geq 3$ and $\gamma_1,\gamma_2\in{\mathbb C}$.
Suppose $\phi:D(Q_2,\gamma_2)\rightarrow D(Q_1,\gamma_1)$ is an isomorphism, and let $f,g,h$ be the images in $D(Q_1,\gamma_1)$ of the standard generators for $D(Q_2,\gamma_2)$.
If $f\in {\mathbb C}[u]$, then $f=u$, $Q_1=Q_2$ and either:

(i) $\gamma_1=\gamma_2$ and $g=v,h=w$ (the trivial isomorphism) or,

(ii) $\gamma_1=-\gamma_2$ and $g=-v$, $h=-w$.
\end{lemma}

\begin{proof}
Since the cosets of $(\ad u)^j(v^m)$, $j\geq 0$ form a basis for $A_m^{(\infty)}/A_{m-1}^{(\infty)}$, the centralizer of $u$ in $D(Q_1,\gamma_1)$ is ${\mathbb C}[u]$.
But therefore the centralizer of $f$ in $D(Q_1,\gamma_1)$ is ${\mathbb C}[f]$.
It follows that $f=au+b$ for some $a\in{\mathbb C}^\times$, $b\in{\mathbb C}$.
By (the proof of) Lemma \ref{Fprod}, $F_1(\ad f,l_f)(h)=0$.
Suppose $g=\xi' u^iv^j+$lower terms (resp. $g=\xi' u^iwv^{j-1}+$lower terms).
Then $h=\xi'' u^{i}wv^{j-1}+$lower terms (resp. $h=\xi'' u^{i+1}v^{j}+$lower terms).
Thus $F_j(\ad u,l_u)(h)\in A_{j-1}^{(\infty)}$ and if $P\in {\mathbb C}[S,T]$ is any polynomial of the form $S^2+c_1S+c_2+dT$ which is not equal to $F_j(\ad u,l_u)$, then $P(\ad u,l_u)(h)\not\in A_{j-1}^{(\infty)}$.
(This is clear since $h,[u,h],uh$ are linearly independent over $A_{j-1}^{(\infty)}$.)
Hence we must have $F_1(\ad f,l_f)=a^2 F_j(\ad u,l_u)$.
It follows that $a=1/j^2$ and $b=1/4-1/4j^2$.
But by the same argument for the inverse isomorphism $\phi^{-1}:D(Q_1,\gamma_1)\rightarrow D(Q_2,\gamma_2)$, we must have $u=f/k^2+(1-1/k^2)$ for some $k\geq 1$.
It follows that $k=j=1$.
Hence $f=u$.

Now $g=\xi' u^iv+$lower terms or $g=\xi' u^iw+$lower terms for some $i\geq 0$ and some $\xi'\in{\mathbb C}^\times$ by Lemma \ref{Fprod}.
Applying the same argument to $\phi^{-1}$, there is an equality $v=q_1(u) g+q_2(u)h+q_3(u)$ for some polynomials $q_1(t),q_2(t),q_3(t)$.
But if $g=\xi' u^iv+$lower terms (resp. $g=\xi' u^iw+$lower terms) then $h=\xi' u^iw+$lower terms (resp. $h=-\xi' u^{i+1}v+$lower terms).
In other words, the leading terms of $q_1(u)g$ and $q_2(u)h$ are of different degrees.
Hence $g=\xi' v+p(u)$ for some polynomial $p(t)$.
Moreover, $F_1(\ad u,l_u)(g)=2\xi'\gamma_1+4up(u)=2\gamma_2$, hence $p(t)=0$.
Thus $g=\xi' v$ and $h=\xi' w$.
But now $$Q_1(u)+uv^2+w^2-2wv-\gamma_1 v=Q_2(u)+(\xi')^2uv^2+(\xi')^2w^2-2(\xi')^2wv-\xi'\gamma_2v=0$$
It follows that $\xi'^2Q_1(u)-\xi'^2\gamma_1 v=Q_2(u)-\xi'\gamma_2 v$, hence that $\xi'=\pm 1$.
This leaves only one non-trivial possibility: that $g=-v$ and $h=-w$.
But one can clearly define such an isomorphism $D(Q_1,-\gamma_1)\rightarrow D(Q_1,\gamma_1)$.
This completes the proof of the Lemma.
\end{proof}

We have therefore solved the isomorphism problem in type $D_{n+1}$, $n\geq 4$.

\begin{definition}
Let $Q(t)$ be a monic polynomial of degree $n\geq 4$, let $\gamma\in{\mathbb C}$ and let $u,v,w$ be the standard generators for $D(Q,\gamma)$.
Then we denote by $\Theta$ the isomorphism $D(Q,\gamma)\rightarrow D(Q,-\gamma)$ which maps $u\mapsto u'$, $v\mapsto -v'$, $w\mapsto -w'$ where $u',v',w'$ are the standard generators for $D(Q,-\gamma)$.
\end{definition}

\begin{rk}
This definition of $\Theta$ should perhaps refer to the defining parameters $Q,\gamma$ in its definition.
However, $\Theta$ can be thought of as the action of the non-identity element of the normalizer $N_{\SL(V)}(\Gamma)$ on the space of noncommutative deformations of $V/\Gamma$.
\end{rk}

\begin{theorem}\label{geq4}
Let $Q(t)$ be a monic polynomial of degree $n\geq 4$ and let $\gamma\in{\mathbb C}$.

(a) If $\tilde{Q}$ is monic of degree greater than or equal to 3 and $\tilde{\gamma}\in{\mathbb C}$, then $D(Q,\gamma)$ is isomorphic to $D(\tilde{Q},\tilde{\gamma})$ if and only if $\tilde{Q}=Q$ and $\tilde\gamma=\pm\gamma$.

(b) The automorphism group of $D(Q,\gamma)$ is trivial unless $\gamma=0$, in which case the automorphism group is cyclic of order 2, generated by $\Theta$.
\end{theorem}

\begin{proof}
This follows from Lemma \ref{outer}, Lemma \ref{Ngeq4}, Cor. \ref{v} and Cor. \ref{iiv}.
\end{proof}

\begin{corollary}\label{moduli}
The moduli space of isomorphism classes of noncommutative deformations of a Kleinian singularity of type $D_n$, $n\geq 5$ is isomorphic to a vector space of dimension $n$.
\end{corollary}

\begin{proof}
The vector space $V$ of monic polynomials of degree $(n-1)$ is isomorphic to ${\mathbb C}^{n-1}$.
Hence we map the isomorphism class of $D(Q,\gamma)$ to $(Q,\gamma^2)\in V\oplus{\mathbb C}\cong{\mathbb C}^n$.
\end{proof}

We apply this to determine when to of the algebras $H(P,\gamma)$ ($P(t)$ has leading term $nt^{n-1}$, $n\geq 4$) are isomorphic.

\begin{theorem}\label{Hn4}
Let $P(t)$ be a polynomial with leading term $nt^{n-1}$ ($n\geq 4$), $\tilde{P}(t)$ a polynomial with leading term $Nt^{N-1}$ ($N\geq 3$) and let $\gamma,\tilde\gamma\in{\mathbb C}$.
Then $H(P,\gamma)\cong H(\tilde{P},\tilde\gamma)$ if and only if $P=\tilde{P}$ and $\gamma=\pm\tilde\gamma$.
\end{theorem}

\begin{proof}
Suppose there exists some isomorphism $\phi:H(P,\gamma)\rightarrow H(\tilde{P},\tilde\gamma)$.
Let $Q(t)$ (resp. $\tilde{Q}(t)$) be the unique monic polynomial with zero constant term such that $Q(-s(s+1))+(s+1)P(-s(s+1))$ (resp. $\tilde{Q}(-s(s+1))+(s+1)\tilde{P}(-s(s+1))$) is an even polynomial in $s$.
Let $\Omega=Q(U)+UV^2+W^2-2WV-\gamma V$ (resp. $\tilde\Omega=\tilde{Q}(\tilde{U})+\tilde{U}\tilde{V}^2+\tilde{W}^2-2\tilde{W}\tilde{V}-\tilde\gamma\tilde{V}$), where $U,V,W$ (resp. $\tilde{U},\tilde{V},\tilde{W}$) are the standard generators for $H(P,\gamma)$ (resp. $H(\tilde{P},\tilde\gamma)$).
By Lemma \ref{centreisp}, $Z(H(P,\gamma))={\mathbb C}[\Omega]$ and $Z(H(\tilde{P},\tilde\gamma))={\mathbb C}[\tilde\Omega]$.
But therefore $\phi(\Omega)=a\tilde\Omega+c$ for some $a\in{\mathbb C}^\times$, $c\in{\mathbb C}$.
It follows that $\phi$ induces an isomorphism $H(P,\gamma)/(\Omega)\rightarrow H(\tilde{P},\tilde\gamma)/(\tilde\Omega-c/a)$.
But $H(P,\gamma)/(\Omega)\cong D(Q,\gamma)$ and $H(\tilde{P},\tilde\gamma)/(\tilde\Omega-c/a)\cong D(\tilde{Q}-c/a,\tilde\gamma)$.
It follows that $\tilde\gamma=\pm\gamma$ and $\tilde{Q}=Q+c/a$, hence $\tilde{P}=P$.
\end{proof}

\section{Isomorphisms in type $D_4$}

Thm. \ref{geq4} solves the problem of determining all isomorphisms $D(Q_2,\gamma_2)\cong D(Q_1,\gamma_1)$ where the degree of $Q_2$ is greater than or equal to 4.
Hence we have only to deal with the case $N=3$.
On considering the inverse isomorphism, we see that $n=3$ as well.
Furthermore, if $\phi:D(Q_2,\gamma_2)\rightarrow D(Q_1,\gamma_1)$ is not of the form described in Lemma \ref{outer} then by Lemma \ref{Ngeq4}, $f=\xi v+p(u)$, $g=\xi' v+q(u)$, $h=\xi'' w+$lower terms and $u=c_1 g+c_2f+c_3$ for some polynomials $p(t),q(t)$ and some $c_1,c_2\in{\mathbb C}^\times,c_3\in{\mathbb C}$.
Replacing $\phi$ by $\phi^{-1}$, we may assume that $p$ is linear.
(We will see that this holds for both $\phi$ and $\phi^{-1}$.)
Since $Q_2(f)+fg^2+h^2-2hg-\gamma_2 g=0$, we must have $\xi'=\pm i\xi$.
After composing with an isomorphism of the form given in Lemma \ref{outer}, if necessary, we may assume furthermore that $\xi'=i\xi$.
Now $2\xi'' w=2h=[f,g]=i\xi [p(u),v]-\xi[q(u),v]$.
It follows that $q$ is also linear.
Hence assume $p(t)=at+b$ and $q(t)=ct+d$.
Then $\xi''=\xi(ia-c)$.
Moreover, $$\begin{array}{rl}[f,h]= & \xi^2(ia-c)[v,w]+\xi a(ia-c)[u,w] \\ = & \xi^2 (ia-c)(v^2+P_1(u))+a(ia-c)\xi(-2uv+2w+\gamma_1)\end{array}$$ where $P_1(t)$ is the unique polynomial such that $Q_1(-s(s+1))+(s+1)P_1(-s(s+1))$ is even in $s$.
On the other hand, by assumption $[f,h]=-2fg+2h+\gamma_2$ and $fg=i\xi^2 v^2+((ai+c)u+(bi+d))\xi v-2ia\xi w+(au+b)(cu+d)$.
We deduce that $a=-1/2$, $c=3i/2$, $d=-bi$ and $-4i\xi^2 P_1(u)+2i\xi\gamma_1=i(u-2b)(3u-2b)+2\gamma_2$.
Suppose $P_1(t)=3t^2+X_1u+Y_1$.
Then it follows that $\xi^2=-1/4$, $8b=-X_1$ and $\gamma_2=i(Y_1/2-X_1^2/32)+i\xi\gamma_1$.
We choose $\xi=i/2$.
(The case $\xi=-i/2$ will then arise as the inverse of the isomorphism we construct below, composed with the non-trivial isomorphism from Lemma \ref{outer}.)
Assume therefore that:
$$f=iv/2-u/2-X_1/8,\; g=-v/2+3iu/2+iX_1/8,\; h=w$$
We wish to determine for which values of $Q_2,\gamma_2$ there exists an isomorphism $\phi$ mapping the standard generators for $D(Q_2,\gamma_2)$ onto $f,g,h$.
By the calculation above, we must have $\gamma_2=i(Y_1/2-X_1^2/32)-\gamma_1/2$.
We note the following description of the coefficients of the polynomial $P(t)$ in terms of those of $Q(t)$.

\begin{lemma}\label{constants}
Suppose $Q(t)=t^3+At^2+Bt+C$ and $P(t)=3t^2+Xt+Y$.
Then $Q(-s(s+1))+(s+1)P(-s(s+1))$ is an even polynomial in $s$ if and only if $X=2A+8$ and $Y=2A+B+8$.
\end{lemma}

\begin{proof}
We have $$Q(-s(s+1))=-s^6-3s^5+(A-3)s^4+(2A-1)s^3+(A-B)s^2-Bs+C\mbox{, and}$$ $$(s+1)P(-s(s+1))=3s^5+9s^4+(9-X)s^3+(3-2X)s^2+(Y-X)s+Y$$
The Lemma follows on comparing odd powers of $s$.
\end{proof}

To proceed, we therefore calculate $Q_2(f)+g^2f+h^2+2gh-\gamma_2 g$, assuming $Q_2(f)=f^3+A_2 f^2+B_2f+C_2$.
By a straightforward calculation \begin{equation}\label{fsquared}f^2=-v^2/4-(u+X_1/4)iv/2+iw/2+(u+X_1/4)^2/4\end{equation} and \begin{equation}\label{gsquared}g^2=v^2/4-(3u+X_1/4)iv/2+3iw/2-(3u+X_1/4)^2/4\end{equation}
Adding (\ref{fsquared}) and (\ref{gsquared}), we obtain $f^2+g^2=-2(u+X_1/8)iv+2iw-2u(u+X_1/8)$.
Multiplying on the right by $f$, we obtain $f^3+g^2f=(u+X_1/8)v^2-wv+((X_1/4-2)u+X_1^2/32)iv-(3u+(X_1/2-2))iw+u(u+X_1/8)(u+X_1/4)+i\gamma_1$.
Hence $f^3+g^2f+h^2=wv+(X_1/8)v^2+((X_1/4-2)u+ X_1^2/32-i\gamma_1)iv-(3u+(X_1/2-2))iw+u(u+X_1/8)(u+X_1/4)-Q_1(u)+i\gamma_1$.
Moreover, $X_1/8=1+A_1/4$.
We deduce that \begin{eqnarray}\nonumber\lefteqn{f^3+g^2f+h^2+2gh=A_1v^2/4+(A_1u/2+2(1+A_1/4)^2-i\gamma_1)iv} \\ \label{long}& & -A_1iw/2+u(u+1+A_1/4)(u+2+A_1/2)-Q_1(u)-P_1(u)+i\gamma_1\end{eqnarray}
It follows that $A_2=A_1$.
Multiplying (1) by $A_1$ and substituting for $X_1$, we have:
\begin{eqnarray}\nonumber\lefteqn{A_1f^2=-A_1v^2/4-(A_1u/2+A_1(1+A_1/4))iv} \\ \label{A1f2} & & +A_1iw/2+A_1u^2/4+A_1(1+A_1/4)u+A_1(1+A_1/4)^2\end{eqnarray}
We notice moreover that $\gamma_2=i(B_1/2+2(1-A_1^2/16)+i\gamma_1/2)$.
It follows that
\begin{eqnarray}\label{gamma2g}-\gamma_2 g=(B_1/4+1-A_1^2/16+i\gamma_1/4)iv+(B_1/4+1-A_1^2/16+i\gamma_1/4)(3u+2+A_1/2)\end{eqnarray}
Let $B_2=6(A_1^2/16-1)+3i\gamma_1/2-B_1/2$.
Taking the sum of (\ref{long}), (\ref{A1f2}) and (\ref{gamma2g}), we see that
\begin{eqnarray}\nonumber\lefteqn{f^3+A_1f^2+g^2f+h^2+2gh-\gamma_2 g=-B_2 iv/2} \\ \label{b1}& \mbox{} & +B_2u/2+B_2(1+A_1/4)+(B_1-i\gamma_1+4(1-A_1^2/16))A_1/4-C_1\end{eqnarray}
Adding $B_2 f$ to (\ref{b1}), we deduce that $Q_2(t)=t^3+A_1 t^2+B_2 t+C_2$ where $C_2=C_1-A_1(B_1/4-i\gamma_1/4+1-A_1^2/16)$.

\begin{lemma}
Let $f=iv/2-u/2-(1+A_1/4)$, $g=-v/2+3iu/2+i(1+A_1/4)$ and $h=w$.
Then $f^3+A_1f^2+B_2 f+C_2+fg^2+h^2-2hg-\gamma_2 g=0$, where $$B_2=6(A_1^2/16-1)+3i\gamma_1/2-B_1/2,$$ $$C_2=C_1-A_1(B_1/4-i\gamma_1/4+1-A_1^2/16),$$ $$\gamma_2=iB_1/2-2i(A_1^2/16-1)-\gamma_1/2.$$
Moreover, $[f,g]=2h$, $[f,h]=-2fg+2h+\gamma_2$ and $[g,h]=g^2+3f^2+(2A_1+8)f+(2A_1+B_2+8)$.
\end{lemma}

\begin{proof}
Let $Q_2(t)=t^3+A_1t^2+B_2t+C_2$.
By construction, $[f,g]=2h$ and $[f,h]=-2fg+2h+\gamma_2$.
Moreover, by the discussion above, $Q_2(f)+g^2 f+h^2+2gh-\gamma_2 g=0$.
But it follows from the commutator relation $[f,g]=2h$ that $g^2f+2gh=fg^2-2hg$.
Hence we have only to show that $[g,h]=g^2+3f^2+(2A+8)f+2A+B+8$.
We deduce from the Jacobi identity $[f,[g,h]]=[[f,g],h]+[g,[f,h]]$ and the known commutator relations for $f$ that $[f,[g,h]]=2(gh+hg)=[f,g^2]$.
Hence $[g,h]=g^2+z$ for some $z\in Z_{D(Q_1,\gamma_1)}(f)=Z(f)$.
We claim that $Z(f)={\mathbb C}[f]$.
Indeed, let $x$ be an element of $Z(f)\setminus{\mathbb C}[f]$ of minimal degree.
By considering $\{ \gr f,\gr x\}$ we see that $x=\chi v^j+$lower terms for some $\chi\in{\mathbb C}^\times$ and $j$.
But now $x-\chi(-2if)^j$ is an element of $Z(f)\setminus{\mathbb C}[f]$ of lower degree than $x$, which provides a contradiction.
It follows that $z=p(f)$ for some polynomial $p(t)$.
Thus there is a homomorphism $H(p,\gamma_2)\rightarrow D(Q_1,\gamma_1)$ which sends the standard generators $U,V,W$ for $H(p,\gamma)$ to $f,g,h$.
The equality $Q_2(f)+fg^2+h^2-2hg-\gamma_2 g=0$ now implies that $p_2(t)=3t^2+(2A_1+8)t+2A_1+B_2+8$.
\end{proof}

\begin{definition}\label{isos}
For a monic polynomial $Q$ of degree $3$ and $\gamma\in{\mathbb C}$, let $\Theta$ be the isomorphism $D(Q,\gamma)\rightarrow D(Q,-\gamma)$ given by $u\mapsto u',v\mapsto -v',w\mapsto -w'$, where $u',v',w'$ are the standard generators of $D(Q,-\gamma)$.

Let $\tilde{B}=6(A^2/16-1)+3i\gamma/2-B/2$, $\tilde{C}=C-A(B/4-i\gamma/4+1-A^2/16)$, $\tilde\gamma=i(2(1-A^2/16)+B/2+i\gamma/2)$ and $\tilde{Q}(t)=t^3+At^2+\tilde{B}t+\tilde{C}$.
Let $\Psi$ be the isomorphism $D(Q,\gamma)\rightarrow D(\tilde{Q},\tilde\gamma)$ given by $$u\mapsto i\tilde{v}/2-\tilde{u}/2-(1+A/4),\; v\mapsto -\tilde{v}/2+3i\tilde{u}/2+i(1+A/4),\; w\mapsto \tilde{w}$$
where $\tilde{u},\tilde{v},\tilde{w}$ are the standard generators of $D(\tilde{Q},\tilde\gamma)$.
\end{definition}

\begin{rk}
As in the case $n\geq 4$, the isomorphisms $\Theta,\Psi$ depend on the choice of $Q,\gamma$ and therefore should perhaps refer to these defining parameters in their definition.
However, one can think of $\Theta$ and $\Psi$ as representatives of elements of $N_{\SL(V)}(\Gamma)/\Gamma$ acting as transformations of the space of noncommutative deformations of $V/\Gamma$.
Since any element of $N_{\SL(V)}$ preserves the invariant ring ${\mathbb C}[V]^\Gamma$, there is a natural action of $N_{\SL(V)}(\Gamma)/\Gamma$ on $V/\Gamma$.
Our construction above and Thm. \ref{main} below therefore say that each such element of $N_{\SL(V)}/\Gamma$ has an induced action on the space of noncommutative deformations of $V/\Gamma$, and this induced action produces all possible isomorphisms between points $(Q,\gamma)$.
\end{rk}

\begin{lemma}
(a) $\Psi^3$ is the identity map on each $D(Q,\gamma)$.

(b) $\Theta\circ\Psi\circ\Theta^{-1}=\Psi^2$.
\end{lemma}

\begin{proof}
Consider $\Psi:D(Q,\gamma)\rightarrow D(\tilde{Q},\tilde\gamma)$ and $\Psi:D(\tilde{Q},\tilde\gamma)\rightarrow D(\hat{Q},\hat\gamma)$.
Let $u,v,w$ (resp, $\tilde{u},\tilde{v},\tilde{w}$, $\hat{u},\hat{v},\hat{w}$) be the standard generators for $D(Q,\gamma)$ (resp. $D(\tilde{Q},\tilde\gamma)$, $D(\hat{Q},\hat\gamma)$).
Then by calculation $\Psi^2:u\mapsto -i\hat{v}/2-\hat{u}/2-1-A/4$, $v\mapsto -\hat{v}/2-3i\hat{u}/2-i(1+A/4)$ and $w\mapsto\hat{w}$.
This proves (b).
On considering the composition of $\Psi^2:D(Q,\gamma)\rightarrow D(\hat{Q},\hat\gamma)$ with $\Psi:D(\hat{Q},\hat\gamma)\rightarrow D(\overline{Q},\overline\gamma)$, we see that $\Psi^3:u\mapsto \overline{u},v\mapsto\overline{v},w\mapsto\overline{w}$.
But therefore $\overline{Q}=Q$ and $\overline\gamma=\gamma$.
\end{proof}

We therefore define the isomorphism $\Psi^{-1}=\Psi^2$: for $A,B,C,\gamma\in{\mathbb C}$ let $\hat{B}=6(A^2/16-1)-3i\gamma/2-B/2$, $\hat{C}=C-A(B/4+i\gamma/4+1-A^2/16)$, $\hat{\gamma}=i(2(A^2/16-1)-B/2)-\gamma/2$ and let $\hat{Q}(t)=t^3+At^2+\hat{B}t+\hat{C}$.
Then there exists an isomorphism $\Psi^{-1}: D(Q,\gamma)\rightarrow D(\hat{Q},\hat\gamma)$ given by $u\mapsto -i\hat{v}/2-\hat{u}/2-(1+A/4)$, $v\mapsto -\hat{v}/2-3i\hat{u}/2-i(1+A/4)$, $w\mapsto\hat{w}$, where $\hat{u},\hat{v},\hat{w}$ are the standard generators for $D(\hat{Q},\hat\gamma)$.

Hence, we have completed our task.

\begin{theorem}\label{main}
Let $Q$ be a monic polynomial of degree $3$ and let $\gamma\in{\mathbb C}$.

(a) There are exactly six isomorphims from $D(Q,\gamma)$ to algebras $D(\overline{Q},\overline\gamma)$, namely $\Id_{D(Q,\gamma)}$ and $\Psi,\Psi^{-1},\Theta,(\Theta\circ\Psi),(\Theta\circ\Psi^{-1})$.

(b) If $\gamma=0$ and $B=4(A^2/16-1)$ then $\Aut D(Q,\gamma)$ is isomorphic to the symmetric group $S_3$, and is generated by $\Psi$ and $\Theta$.

(c) $\Aut D(Q,\gamma)$ is of order 2 if exactly one of $B-4(A^2/16-1)-i\gamma$, $B-4(A^2/16-1)+i\gamma$ and $\gamma$ is zero.
If $B=4(A^2/16-1)+i\gamma$ (resp. $B=4(A^2/16-1)-i\gamma$) and $\gamma\neq 0$ then $\Aut D(Q,\gamma)$ is generated by $\Theta\circ\Psi$ (resp. $\Theta\circ\Psi^{-1}$). If $\gamma=0$ but $B\neq 4(A^2/16-1)$ then $\Aut D(Q,\gamma)$ is generated by $\Theta$.

(d) If $\gamma\neq 0$ and $B\neq 4(A^2/16-1)\pm i\gamma$ then there are no non-trivial automorphisms of $D(Q,\gamma)$.
\end{theorem}

\begin{proof}
This follows from Lemma \ref{reducetov}, Cor. \ref{iiv}, Cor. \ref{v}, Lemma \ref{Ngeq4}, Lemma \ref{outer} and the discussion above.
\end{proof}

We therefore think of the symmetric group $S_3$ as acting on the space of noncommutative deformations via the representatives $\Theta,\Psi$.
However, it is not difficult to describe the invariants with respect to this action.
Let $\sigma=\left(\begin{array}{ccc} 1 & 2 & 3 \end{array}\right)$ and $\tau=\left(\begin{array}{cc} 1 & 2 \end{array}\right)$ be generators for $S_3$.

\begin{lemma}\label{invs}
Let $S_3$ act as algebra automorphisms of ${\mathbb C}[A,B,C,\gamma]$ with the action of $\sigma$ (resp. $\tau$) given by that of $\Psi$ (resp. $\Theta$) on $Q(t)=t^3+At^2+Bt+C$ and $\gamma$.
Let $x_1=B-4(A^2/16-1)-\gamma\sqrt{3},x_2=B-4(A^2/16-1)+\gamma\sqrt{3},x_3=6C-AB$ and $x_4=A$.

Then ${\mathbb C}[A,B,C,\gamma]={\mathbb C}[x_1,x_2,x_3,x_4]$ and $\sigma(x_1)=e^{2\pi i/3} x_1,\sigma(x_2)=e^{-2\pi i/3}x_2$, $\tau(x_1)=x_2,\tau(x_2)=x_1$.
The action of $S_3$ on $x_3$ and $x_4$ is trivial.
\end{lemma}

\begin{proof}
The fact that ${\mathbb C}[A,B,C,\gamma]={\mathbb C}[x_1,x_2,x_3,x_4]$ is clear, since $x_2-x_1=2\gamma\sqrt3$, $x_1+x_2+4(x_4^2-1)=2B$ and $x_3+AB=6C$.
The action of $\sigma$ and $\tau$ on $x_1,x_2,x_3,x_4$ follows immediately from definition \ref{isos}.
\end{proof}

\begin{corollary}\label{moduli4}
The moduli space of isomorphism classes of noncommutative deformations of a Kleinian singularity of type $D_4$ is isomorphic to a vector space of dimension 4.
\end{corollary}

\begin{proof}
By Lemma \ref{invs} the ring of invariants ${\mathbb C}[A,B,C,\gamma]^{S_3}$ is generated by $x_1^3+x_2^3,x_1x_2,x_3$ and $x_4$.
But hence the map $D(Q,\gamma)\mapsto ((B-4(A^2/16-1))((B-4(A^2/16-1))^2+9\gamma^2),(B-4(A^2/16-1))^2-3\gamma^2,6C-AB,A)$ induces a bijective map on isomorphism classes of deformations.
\end{proof}

Finally, we can now solve the problem of when two algebras $H(P,\gamma),H(\tilde{P},\tilde\gamma)$ are isomorphic.

\begin{theorem}\label{H4}
Let $P(t)=3t^2+Xt+Y$, $\tilde{P}(t)=3t^2+\tilde{X}t+\tilde{Y}$ and $\gamma,\tilde\gamma\in{\mathbb C}$.
Then $H(P,\gamma)\cong H(\tilde{P},\tilde\gamma)$ if and only if $X=\tilde{X}$ and either

(i) $\tilde{Y}=3(X+X^2/32+i\gamma/2)-Y/2$ and $\pm\tilde\gamma=i(Y/2-X-X^2/32)-\gamma/2$, or

(ii) $\tilde{Y}=3(X+X^2/32-i\gamma/2)-Y/2$ and $\pm\tilde\gamma=-i(Y/2-X-X^2/32)-\gamma/2$, or

(iii) $\tilde{Y}=Y$ and $\tilde\gamma=\pm\gamma$.
\end{theorem}

\begin{proof}
Let $\phi:H(P,\gamma)\rightarrow H(\tilde{P},\tilde\gamma)$ be an isomorphism.
Let $Q(t)$ (resp. $\tilde{Q}(t)$) be the unique monic polynomial with zero constant term such that $Q(-s(s+1))+(s+1)P(-s(s+1))$ (resp. $\tilde{Q}(-s(s+1))+(s+1)\tilde{P}(-s(s+1))$) is even in $s$ and let $\Omega=Q(u)+uv^2+w^2-2wv-\gamma v$, $\tilde\Omega=\tilde{Q}(\tilde{u})+\tilde{u}\tilde{v}^2+\tilde{w}^2-2\tilde{w}\tilde{v}-\tilde\gamma \tilde{v}$.
By Lemma \ref{centreisp}, $Z(H(P,\gamma))={\mathbb C}[\Omega]$ and $Z(H(\tilde{P},\tilde\gamma))={\mathbb C}[\tilde\Omega]$.
It follows that $\phi(\Omega)=a\tilde\Omega+c$ for some $a\in{\mathbb C}^\times$, $c\in{\mathbb C}$.
Hence $\phi$ induces an isomorphism $H(P,\gamma)/(\Omega)\rightarrow H(\tilde{P},\tilde\gamma)/(\tilde\Omega-c/a)$.
But $H(P,\gamma)/(\Omega)\cong D(Q,\gamma)$ and $H(\tilde{P},\tilde\gamma)/(\tilde\Omega-c/a)\cong D(\tilde{Q}-c/a,\tilde\gamma)$.
The theorem now follows from Thm. \ref{main} and the fact that $Q(t)=t^3+(X/2-4)t^2+(Y-X)t$, $\tilde{Q}(t)=t^3+(\tilde{X}/2-4)t^2+(\tilde{Y}-\tilde{X})t$.
\end{proof}

\end{document}